\documentclass[12pt]{article}

\usepackage{amsmath, amsthm, amsfonts, amscd, amssymb}
\usepackage{xcolor}
\usepackage{graphics}
\usepackage{bbm}
\usepackage{dsfont}
\usepackage{verbatim}

\newcommand{\be}{\begin{equation}}
\newcommand{\ee}{\end{equation}}

\DeclareMathOperator*{\argmin}{arg\,min}
\newtheorem{thm}{Theorem}[section]
\usepackage{courier}
\title{Information Content in Data Sets for a Nucleated-Polymerization Model}
\author{H.T. Banks$^1$, M. Doumic$^{2,3}$, C. Kruse$^{2,3}$, S. Prigent$^{2,3}$,  H.Rezaei$^{4}$ \\
{\footnotesize $^1$ Center for Research in Scientific Computation}\\
{\footnotesize North Carolina State University, Raleigh, NC 27695-8212}\\
{\footnotesize $^2$Institut National de Recherche en Informatique et Automatique}\\
 {\footnotesize Paris-Rocquencourt, France}\\
   {\footnotesize $^3$Pierre et Marie Curie University}\\
   {\footnotesize Paris, France}\\
   {\footnotesize $^{4}$Institut National de Recherche Agronomique}\\
    {\footnotesize Jouy-en-Josas, France}\\
 }

\begin{document}
\pagestyle{plain}
\maketitle
\date{}

\begin{abstract}
We illustrate the use of tools (asymptotic theories of standard error quantification using appropriate statistical models, bootstrapping, model comparison techniques) in addition to sensitivity  that may be employed to determine the information content in data sets. We do this in the context of recent models \cite{DouPri2012} for nucleated polymerization in proteins, about which very little is known regarding the underlying mechanisms; thus the methodology we develop here may be of great help to experimentalists.
\end{abstract}
\vspace{.2in}

\vspace{0.2in}

\clearpage

\section{Introduction}


As mathematical models become more complex with multiple states and many parameters to be estimated using experimental data, there is a
need for critical analysis in model validation related to the reliability of parameter estimates obtained in model fitting. A recent concrete example involves
previous HIV models \cite{ABDR,banks_modelling_2008} with 15 or more parameters to be  estimated. In \cite{HIVUQ}, using recently developed parameter selectivity tools \cite{parameterselection} based on parameter sensitivity based scores, it was shown that many  parameters could not be estimated with any degree of reliability. Moreover, we found that quantifiable uncertainty varies among patients depending upon the number of treatment interruptions (perturbations of therapy). This leads to a {\em fundamental question}: how much information with respect to model validation can be expected in a given data set or collection of data sets?

Here we illustrate the use of other tools (asymptotic theories of standard error quantification using appropriate statistical models, bootstrapping, and model comparison techniques) in addition to sensitivity theory that may be used to determine the information content in data sets. We do this in the context of recent models \cite{DouPri2012} for nucleated polymerization in proteins.

After presenting the biological context of amyloid formation, we describe the model in Section~\ref{sec:model}. In Section~\ref{sec:IP}, we investigate the statistical model to be used with our noisy data. This is a necessary step in order to use the correct error model in our  generalized least squares (GLS) minimization. This also reveals information on our experimental observation process. Once we have found parameters which allow a reasonable fit, we determine  the confidence we may have in our estimation procedures. We do this in Section~\ref{sec:SE}, using both the condition number of the covariance matrix and a sensitivity analysis. This reveals a smaller number of parameters (than those estimated  in \cite{DouPri2012}) which appear as reasonably sensitive to the data sets, whereas other do not really affect the quality of the fits to our data. To further support our sensitivity findings, we then apply a bootstrapping analysis in Section~\ref{sec:boot}. We are lead to four main parameters and compare their resulting  errors with the asymptotic confidence intervals of Section~\ref{sec:SE}.
Finally, in Section~\ref{sec:chi2}, we carry out model comparison tests \cite{BanksFitz,BHT2014,BT2009} as used in \cite{MComp}, and these lead us to select three well-defined parameters that can be reliability estimated out of the nine original ones estimated in \cite{DouPri2012}.

\subsection{Protein Polymerization}
It is now known that several neuro-degenerative disorders, including Alzheimer’s disease, Huntington’s disease and Prion diseases e.g., mad cow, are related to aggregations of proteins presenting an abnormal folding. These protein aggregates are called {\em amyloids} and have become a focus of modeling efforts in recent years \cite{CL2010,DouPri2012, SHR2007, Xue2009,Xue2013}.  One of the main challenges in this field is to understand the key  aggregation mechanisms, both qualitatively and quantitatively. In order to test our methodology on a relatively simple case,
 we focus here on  polyglutamine (PolyQ) containing proteins. This was also the case study chosen to illustrate the fairly general ODE-PDE model proposed  in \cite{DouPri2012}; the reason for our choice is that, as shown in \cite{DouPri2012}, the polymerization mechanisms prove to be simpler for PolyQ aggregation than for other types of proteins, e.g. PrP~\cite{Rezaei2007}. To understand data sets from experiments carried by Human Rezaei and his team at INRA, (Virologie et Immunologie Moleculaires), see \cite{DouPri2012}, we adapt the general model to this context.  The data sets (DS1-DS4) of interest to us here are depicted in Figure \ref{datasets} below.

\begin{figure}[h!]
 \begin{center}
 \resizebox{10.5cm}{!}{\includegraphics{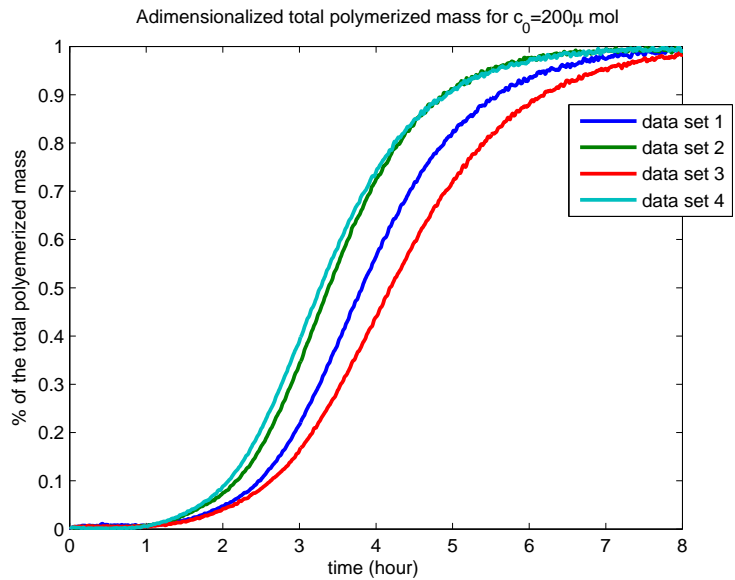}}
 \end{center}
 \caption{The data sets of interest from \cite{DouPri2012,BDK1}.}
 \label{datasets}
 \end{figure}

In \cite{DouPri2012} and a subsequent effort in \cite{BDK1}, the authors sought to investigate several questions including (i) understanding the key polymerization mechanisms, (ii) how to select parameters and calibrate the model, and (iii) how to numerically approximate the model. Here we briefly summarize results related to (iii) and focus primarily on (ii).

\section{The Model}\label{sec:model}
\subsection{Original ODE Model}
This model we used is the same as that of \cite{DouPri2012}. We briefly outline that model. Let $(V, V^*, c_i)$  be the concentrations of the normal monomeric proteins that we will call monomers, of the monomeric proteins presenting an abnormal configuration that we will call conformers, and of the $i$-polymers made of $i$ aggregated abnormal proteins, respectively. The following comprise the fundamental dynamics modeled in \cite{DouPri2012}:
	
	\begin{itemize}
	\item Monomer-conformer exchange: $V \stackrel {k_I^+}{\underset {k_I^-}{\rightleftharpoons}} V^*$
  \item Nucleation: \hspace{0cm} $\underbrace{V^*+V^*+...+V^*}_{i_0}  \stackrel {k_{on}^N}{\underset {k_{off}^N}{\rightleftharpoons}} c_{i_0}$
	\item Polymerization by conformer addition: $c_i+V^* \stackrel {k_{on}^i}{\underset {}{\rightleftharpoons}} c_{i+1}$
	\end{itemize}	
	\vspace{0.2cm}
	Other reactions like fragmentation and coalescence are  negligible for the case of polyglutamine containing proteins (see~\cite{DouPri2012} for experimental justification).

	\vspace{0.5cm}
The law of mass action in the deterministic framework (see \cite{BT2009,Rubinow} and the numerous references therein), translates
	$A+B \stackrel {k_I^+}{\underset {k_I^-}{\rightleftharpoons}} A'+B'$ into the ordinary differential equation
		$\frac{d[A]}{dt}=-k^+ [A] [B] + k^- [A]' [B]'$.

Using these basic ideas we obtain the {\em infinite} system of ordinary differential equations (ODEs) studied in \cite{DouPri2012}
\begin{align}
\frac{dV}{dt}&=-k_I^+ V + k_I^-V^*,\\
\frac{dV^*}{dt}&=k_I^+V-k_I^{-}V^*+i_0 k_{off}^N c_{i_0}-V^*\sum_{i\geq i_0} k_{on}^i c_i,\\
\frac{d c_{i_0}}{dt}&=k_{on}^N(V^*)^{i_0}-k_{off}^N c_{i_0} - k_{on}^{i_0} c_{i_0} V^*,\\
\frac{d c_i}{dt}&=V^*(k_{on}^{i-1}c_{i-1}-k_{on}^i c_i),\hspace{0.5cm} i=i_0+1,.... \label{odeinfinite}
\end{align}	
with initial conditions
\begin{align*}
V(0)=c_0, \hspace{0.2cm} V^*(0)=0,\hspace{0.2cm} c_{i_0}(0)=c_i(0)=0
\end{align*}
and the mass balance equation
\begin{align*}
\frac{d}{dt}\left(V+ V^* + \sum_{i=i_0}^{\infty} i c_i \right)=0.
\end{align*}

The experiments of interest to us measure the total polymerized mass, i.e.,
\begin{align*}
M(t)=\sum_{i\geq i_0} i c_i(t).
\end{align*}

\subsection{An Approximate PDE System and the Associated Forward Problem}
Since very long polymers (a fibril may contain up to $10^6$ monomer units) characterize amyloid formations, a PDE version of the standard model, where a continuous variable $x$ approximates the discrete sizes $i$, is a reasonable approximation for large amyloid polymers. However, for small polymer sizes this curarization does not work very well. Thus we take a "hybrid approach" of leaving the ODE for smaller sizes and use the PDE for larger ones, see \cite{BDK1}.

We define a small parameter $\varepsilon=\frac{1}{i_M}$, and let $x_i=i\varepsilon$ with $i_M \gg 1$ be the average polymer size defined by
\begin{align*}
i_M= \frac{\sum\limits_{i\geq i_0} i c_i }{\sum c_i}.
\end{align*}

Then after definition of dimensionless quantities
\begin{align*}
c^\varepsilon(t,x)=\sum c_i \mathds {1}_{[x_i, x_{i+1}]}
\end{align*}
we may obtain a partial differential equation (PDE) to replace the infinite ODE system.  Rigorous derivations of such continuous integro-PDE models may be found in~\cite{LM2002} for coagulation-fragmentation equations, in~\cite{CGPV02} for the limit of the Becker-D\"oring system toward Lifshitz-Slyozov model, and in~\cite{DGL2009} for the growth-fragmentation "Prion Model". A formal derivation for a full model, also including nucleation, is carried out in~\cite{DouPri2012}.

Let $N_0\in \mathbb{N}$. We then use the approximation
\begin{align}\nonumber
\frac{dV}{dt}&=-k_I^+ V + k_I^-V^*,\\
\frac{dV^*}{dt}&=k_I^+V-k_I^{-}V^*+i_0 k_{off}^N c_{i_0}-V^*\sum_{i\geq i_0} k_{on}^i c_i,\\
\frac{d c_{i_0}}{dt}&=k_{on}^N(V^*)^{i_0}-k_{off}^N c_{i_0} - k_{on}^{i_0} c_{i_0} V^*,\\
{\frac{d c_i}{dt}}&{=V^*(k_{on}^{i-1}c_{i-1}-k_{on}^i c_i),\hspace{0.5cm} i\leq N_0,} \\
{\partial_t c^{\varepsilon}(x,t)}&{=-V^* \partial_x (k_{on} c^{\varepsilon}(x,t)), \hspace{0.5cm} x\geq N_0,}
\end{align}	
with initial conditions
\begin{align*}
V(0)=c_0, \hspace{0.2cm} V^*(0)=0,\hspace{0.2cm} c_{i_0}(0)=c_i(0)=0, \hspace{0.2cm} c^{\varepsilon}(x,0)=0,
\end{align*}
and the boundary condition
\begin{align*}
c^{\epsilon}(x=N_0,t)=c_{N_0}(t).
\end{align*}


Then an assumed mass balance equation becomes

\begin{align*}
\frac{d}{dt}\left(V+ V^* + \sum_{i=i_0}^{N_0} i c_i +\int_{N_0}^{\infty} x c^{\varepsilon}(x)\, dx \right)=0.
\end{align*}



In \cite{BDK1} we considered requirements for a good discretization  scheme including (i) it should conserve the total polymerized mass, (ii) it should be fast and most importantly, (iii) it should be accurate.



To ensure the mass conservation, we replace the ODE for $V^*$ by the {\em mass conservation equation} and obtain
\begin{align}\nonumber
\frac{dV}{dt}&=-k_I^+ V + k_I^-V^*,\\\nonumber
{V^*}&{=c_0-V-\sum_{i=i_0}^{N_0} i c_i - \int_{N_0}^{\infty} x c^{\varepsilon}\,dx} ,\\\nonumber
\frac{d c_{i_0}}{dt}&=k_{on}^N(V^*)^{i_0}-k_{off}^N c_{i_0} - k_{on}^{i_0} c_{i_0} V^*,\\\nonumber
\frac{d c_i}{dt}&=V^*(k_{on}^{i-1}c_{i-1}-k_{on}^i c_i),\hspace{0.5cm} i\leq N_0, \\\nonumber
\partial_t c^{\varepsilon}(x,t)&=-V^* \partial_x (k_{on} c^{\varepsilon}(x,t)), \hspace{0.5cm} x\geq N_0,
\end{align}	
with initial and boundary conditions as before.

We developed methodology for forward solutions in \cite{BDK1}. In considering these forward solutions we first observed that
the desired spatial computational domain is very large as determined by the maximum size of observed polymers, with range up to $10^6$ and the peak in the distribution is at the left side of the domain of interest; for larger polymer sizes, the distribution is almost linearly decreasing.

 Based on these and other considerations discussed in \cite{BDK1}, the PDE was approximated by the Finite Volume Method (see \cite{Lev02} for discussions of Upwind, Lax-Wendroff and flux limiter methods) with an adaptive mesh, refined toward the smaller polymer sizes.
Furthermore, we kept the ratio between the step size and the corresponding mesh element constant, i.e., we used $\frac{\Delta x_i}{x_i}=q<1$ so that   $x_i=\frac{1}{1-q}x_{i-1}$. This mesh is quasi-linear in the sense of  $\frac{\Delta x_{i-1}}{\Delta x_{i}}=1+O(q)$. The resulting Upwind and Lax-Wendroff schemes are then consistent on the progressive mesh (see \cite{Lev02}). For further details on these schemes including examples demonstrating  convergence properties, the interested reader may consult \cite{BDK1}.


\section{The Inverse Problem}
\label{sec:IP}
\begin{figure}[h]
 \begin{center}
 \resizebox{12cm}{!}{\includegraphics{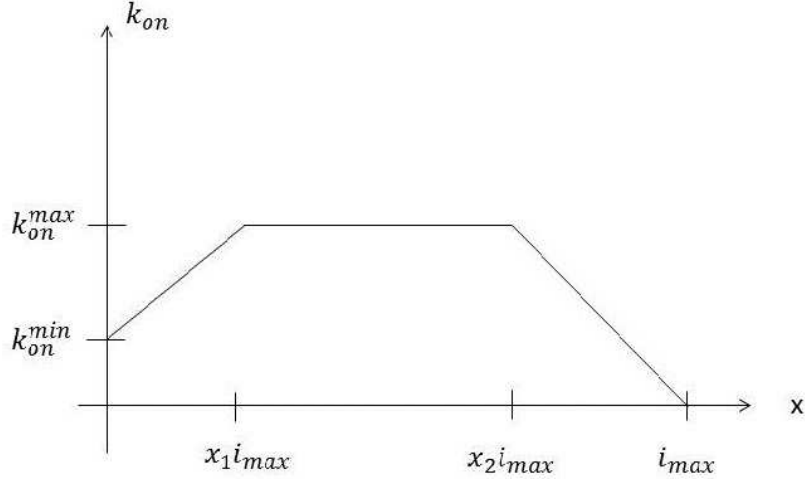}}
 \end{center}
 \caption{Parametric representation for $k_{on}$.\label{fig:kon}}
 \label{kon}
 \end{figure}

A major question in formulating the model for use in inverse problem scenarios consists of how to best parametrically represent the function  $k_{on}$ for our application?  Following~\cite{DouPri2012}, we chose to approximate $k_{on}$ by a function as depicted in Figure \ref{kon}. (According to our discussions between  S. Prigent, H. Rezaei and J. Torrent, other choices like a Gaussian bell curve are also possible, and we  discuss this later). Thus with this parametrization we have 5 more parameters $k_{on}^{min}, k_{on}^{max}, x_1, x_2, i_{max}$ in addition to the 4 basic parameters  $ k_I^+, k_I^-, k_{on}^{N}, k_{off}^{N}$ to be estimated using our data sets.

Thus we seek to estimate (with \underline{\em acceptable quantification of uncertainties}) the  nine parameters $k_I^{+}, k_I^{-}, k_{off}^N, k_{on}^N$, and $k_{on}$ (represented in parametrical form depicted above with the 5 additional unknowns
$k_{on}^{min}, k_{on}^{max},$ $x_1, x_2, i_{max}$) that fit the data best! To do this we need an efficient discretization method as discussed above for the forward problem {\em as well as a correct assumption on the measurement errors} in the inverse problem.

\subsection{Estimation of Parameters}

We make some standard statistical assumptions (see \cite{BHT2014,BT2009,DG,SeWi}) underlying our inverse problem formulations.
\begin{itemize}
	\item Assume that there exists a {\em true or nominal set} of parameter $\theta_0=(k_I^-,...,i_{max})$
	\item Let $\mathcal{E}_i$ be iid with $\mathbb{E}(\mathcal{E}_i)=0$ and cov$(\mathcal{E}_i,\mathcal{E}_j)=\sigma^2$. Let $\epsilon_i \in \mathcal{E}_i$.
\end{itemize}
Denote the {estimated parameter for $\theta_0$ as $\hat{\theta}$.
The inverse problem is based on statistical assumptions on the {observation error} in the data.
\vspace{1cm}

If we assume an {\em absolute error data model} then  data points are taken with equal importance. This is represented
by observations
\begin{equation}\label{absolute}
y_i=M(t_i,\theta_0)+\epsilon_i.
\end{equation}

On the other hand, if one assumes some type of {\em relative error data model} then  the error is proportional in some sense to the measured polymerized mass. This can be represented by observations of the form
\begin{equation}\label{relative}
y_i=M(t_i,\theta_0)+M(t_i,\theta_0)^{\gamma}\epsilon_i, \hspace{0.3cm} \gamma\in (0,1].
\end{equation}



Absolute model error formulations dictate we use {\em Ordinary Least Squares (OLS)} inverse problem \cite{BHT2014,BT2009} given by

\begin{equation}\label{OLSform}
\hat{\theta}=\argmin \sum (y_i-M(t_i,\theta))^2
\end{equation}
while for relative error model  one should use inverse problem formulations with {\em Generalized Least Squares (GLS)} cost functional
\begin{equation}\label{WLSform}
\hat{\theta}=\argmin \sum \left(\frac{y_i-M(t_i,\theta)}{M(t_i,\theta)^{\gamma}}\right)^2, \hspace{0.3cm} \gamma\in (0,1].
\end{equation}


\subsubsection{The Residual Plots}

To obtain a correct statistical model, we used {residual plots} (see \cite{BHT2014,BT2009} for more details) with residuals given by
\begin{align*}
r_i=\frac{y_i-M(t_i,\hat{\theta})}{M(t_i,\hat{\theta})^{\gamma}}, \hspace{0.3cm} \gamma \in [0,1]
\end{align*}

To illustrate what we are seeking for our data sets, we first used \underline{simulated} {\em relative error data} (simulated data for $\gamma =1$), then carried out the inverse problems for both a relative error cost functional (i.e., $\gamma=1$) and an ordinary least squares cost functional (i.e., $\gamma=0$). We then plotted the corresponding residuals vs time and also residuals vs the model values. The first plots are related to the correctness of our assumption of independency and identical distributions {\em i.i.d.} for the data whereas the second plots contain information as to the correctness of the form of our proposed statistical model.

 \begin{figure}[h!]
 \begin{tabular}{cc}
 \resizebox{5.6cm}{!}{\includegraphics{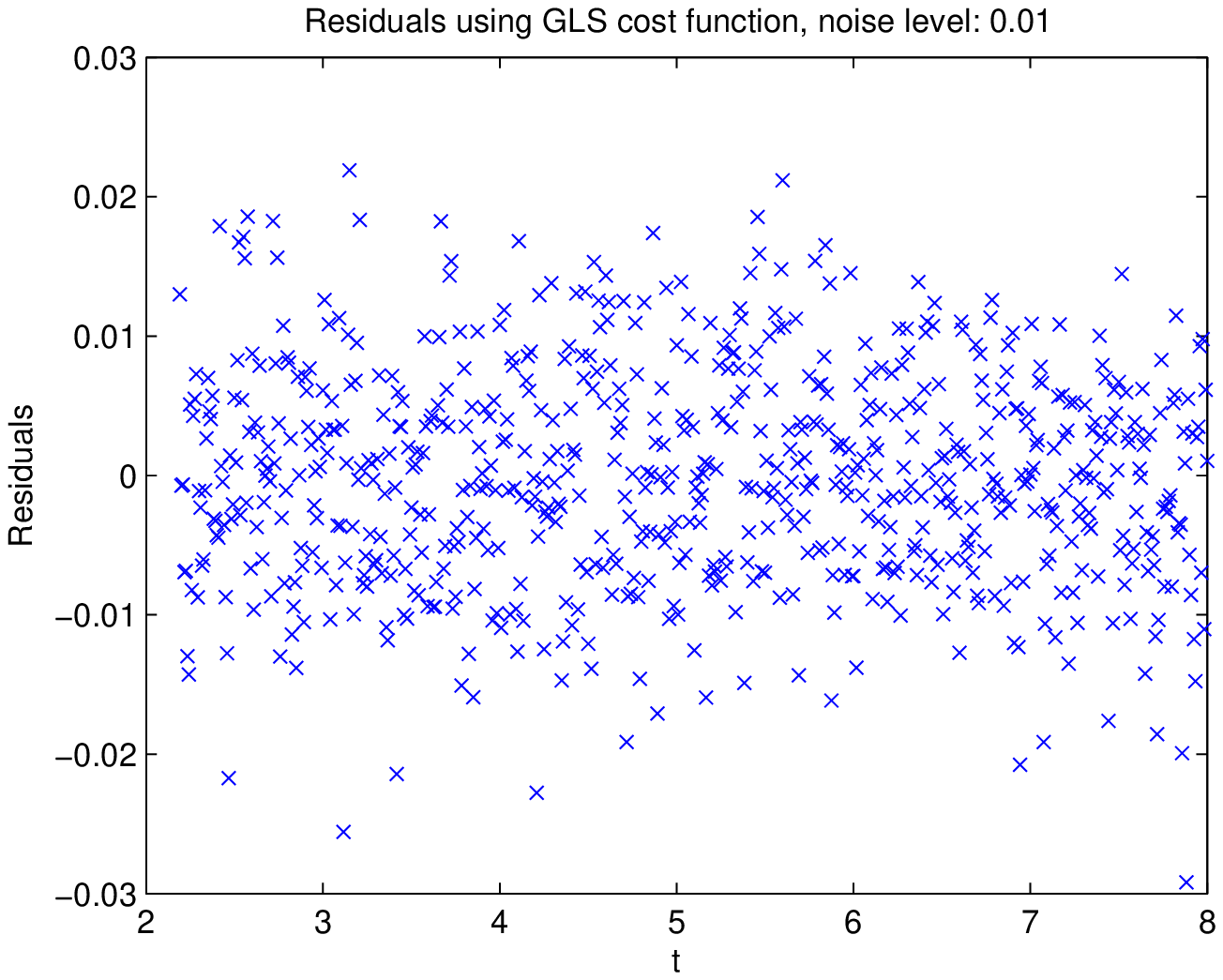}}&
 \resizebox{5.6cm}{!}{\includegraphics{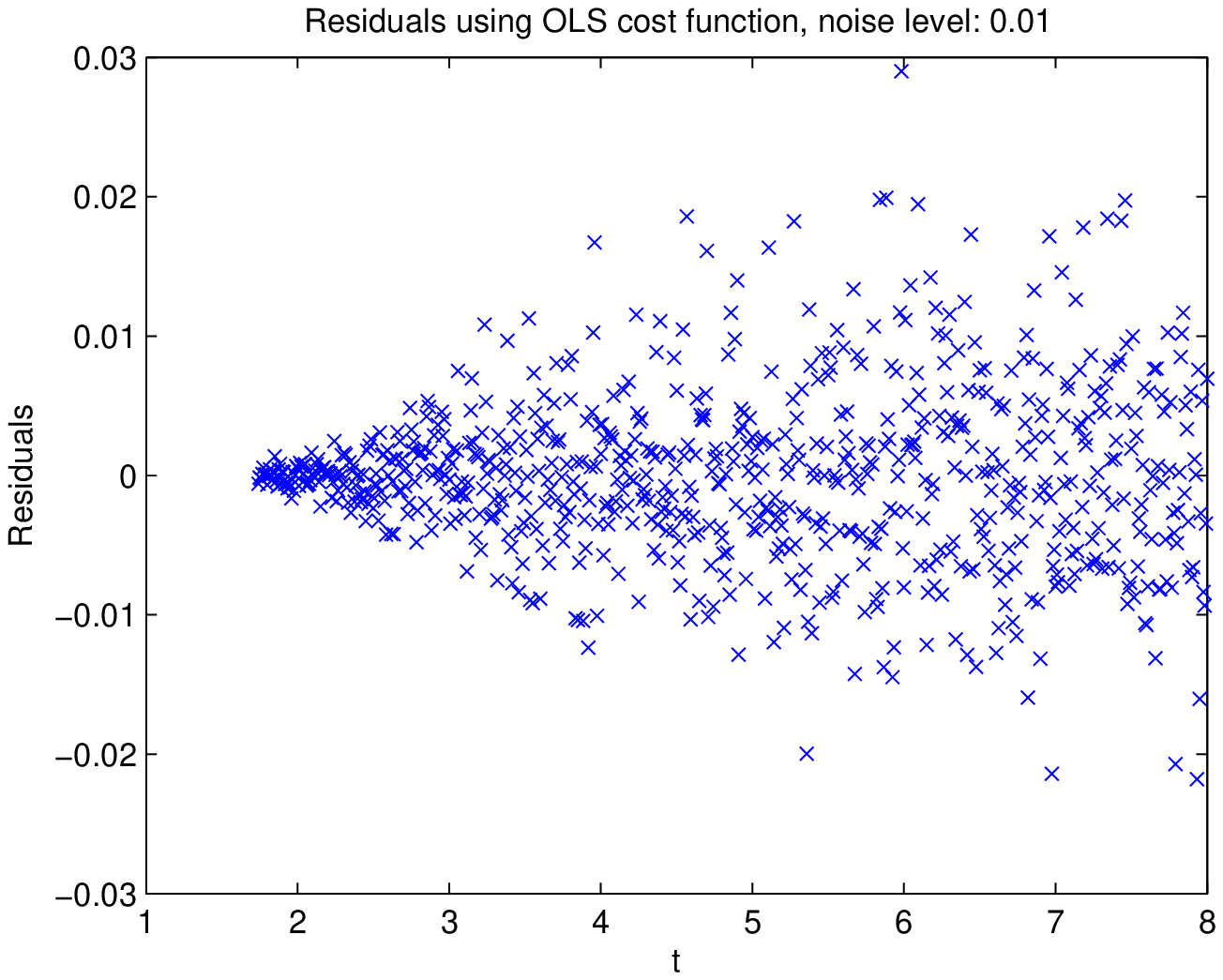}}\\
 (a) & (b)
 \end{tabular}
 \caption{ Plots with simulated data: (a) Correct cost function vs. time $(\gamma=1)$; (b)Incorrect cost function vs. time $(\gamma=0)$ }
 \end{figure}

 \begin{figure}[h!]
 \begin{tabular}{cc}
 \resizebox{5.6cm}{!}{\includegraphics{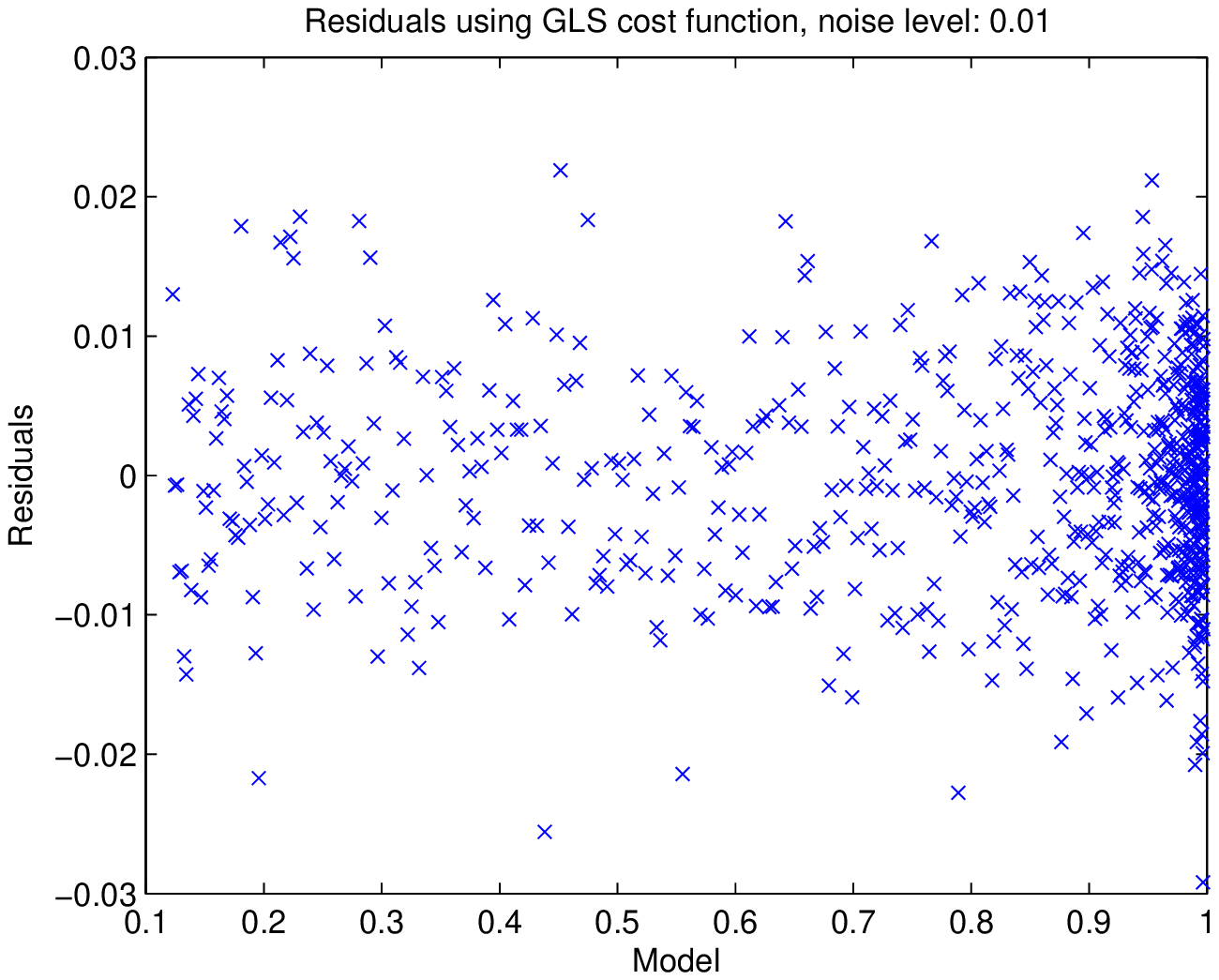}}&
 \resizebox{5.6cm}{!}{\includegraphics{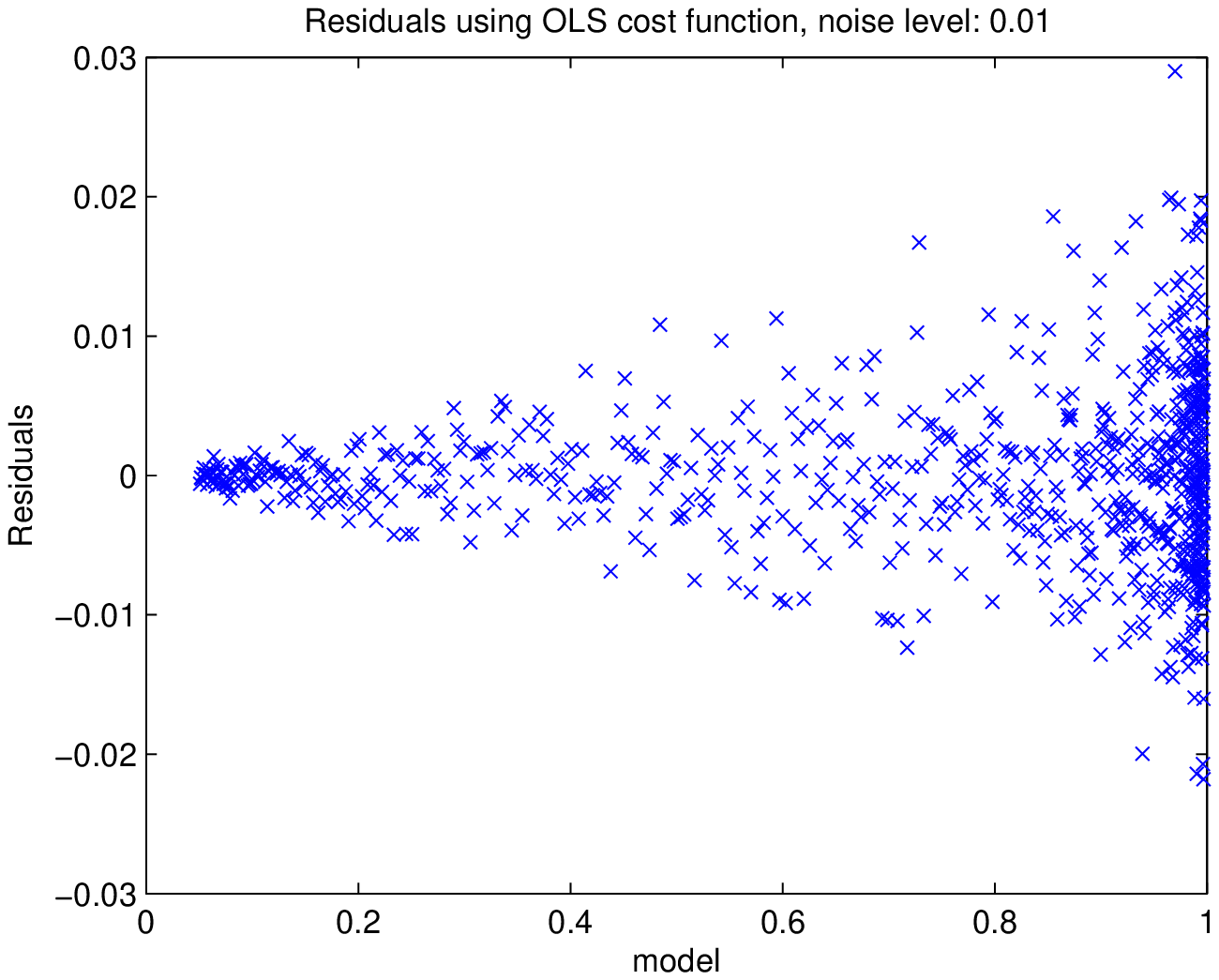}}\\
 (a) & (b)
 \end{tabular}
 \caption{ Plots with simulated data: (a) Correct cost function vs. model $(\gamma=1)$; (b)Incorrect cost function vs. model $(\gamma=0)$ }
 \end{figure}

\clearpage

\subsection{Statistical Models of Noise}
We next carried out similar inverse problems with data set (DS) 4 of our experimental data collection. We first used DS 4 on the interval $t\in [0,8]$. Based on some earlier calculations we also chose the nucleation index $i_0=2$ for all our subsequent calculations. The residual plots given below in Figures \ref{OLSdataset4} and \ref{GLSdataset4} suggest strongly that \underline{\em neither} of the first attempts of assumed statistical models and corresponding cost functionals (absolute error and OLS or relative error with $\gamma=1$ and simple GLS) are correct.
 \begin{figure}[h!]
 \begin{center}
 \begin{tabular}{cc}
\resizebox{6.8cm}{!}{\includegraphics{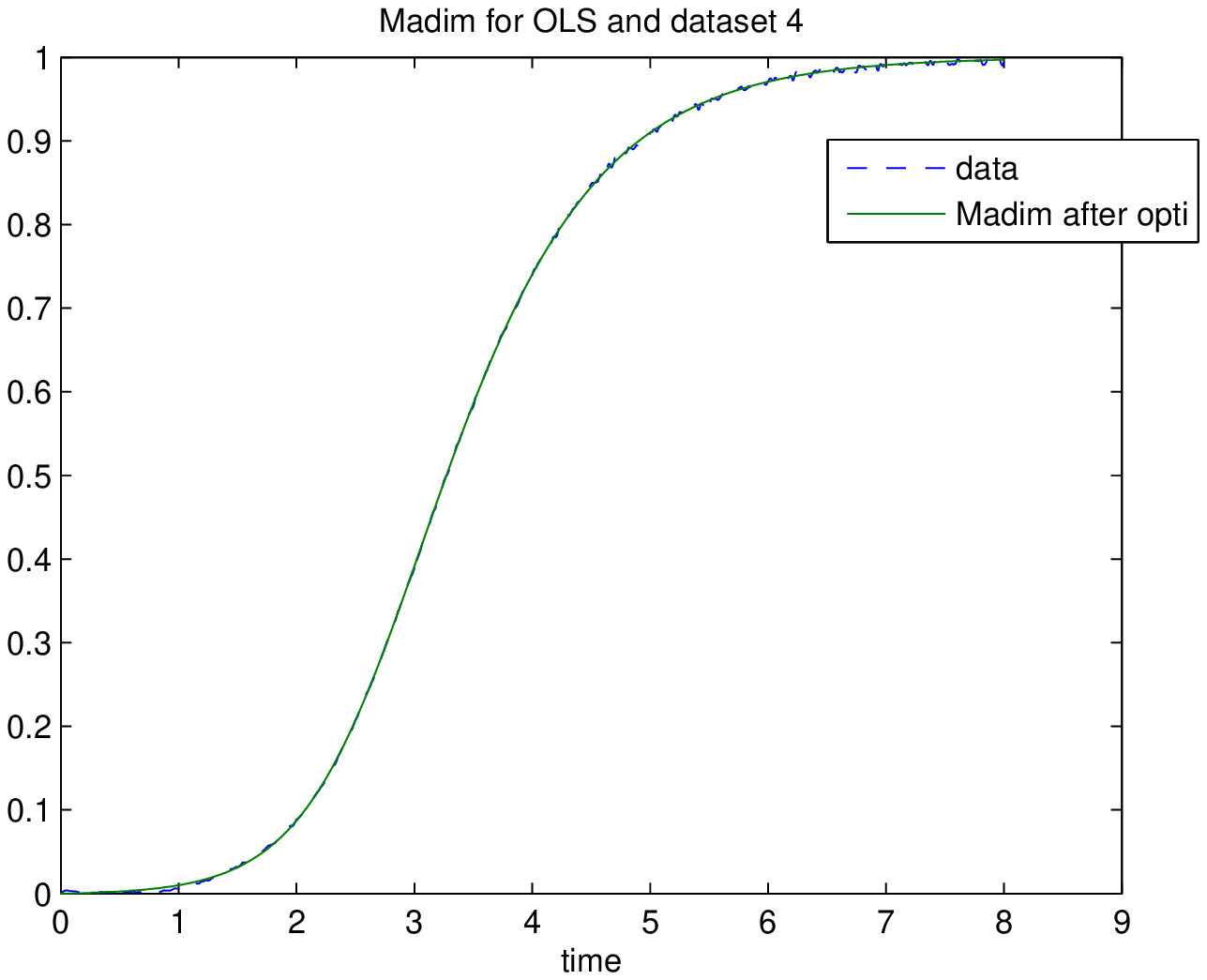}}&
 \resizebox{6.8cm}{!}{\includegraphics{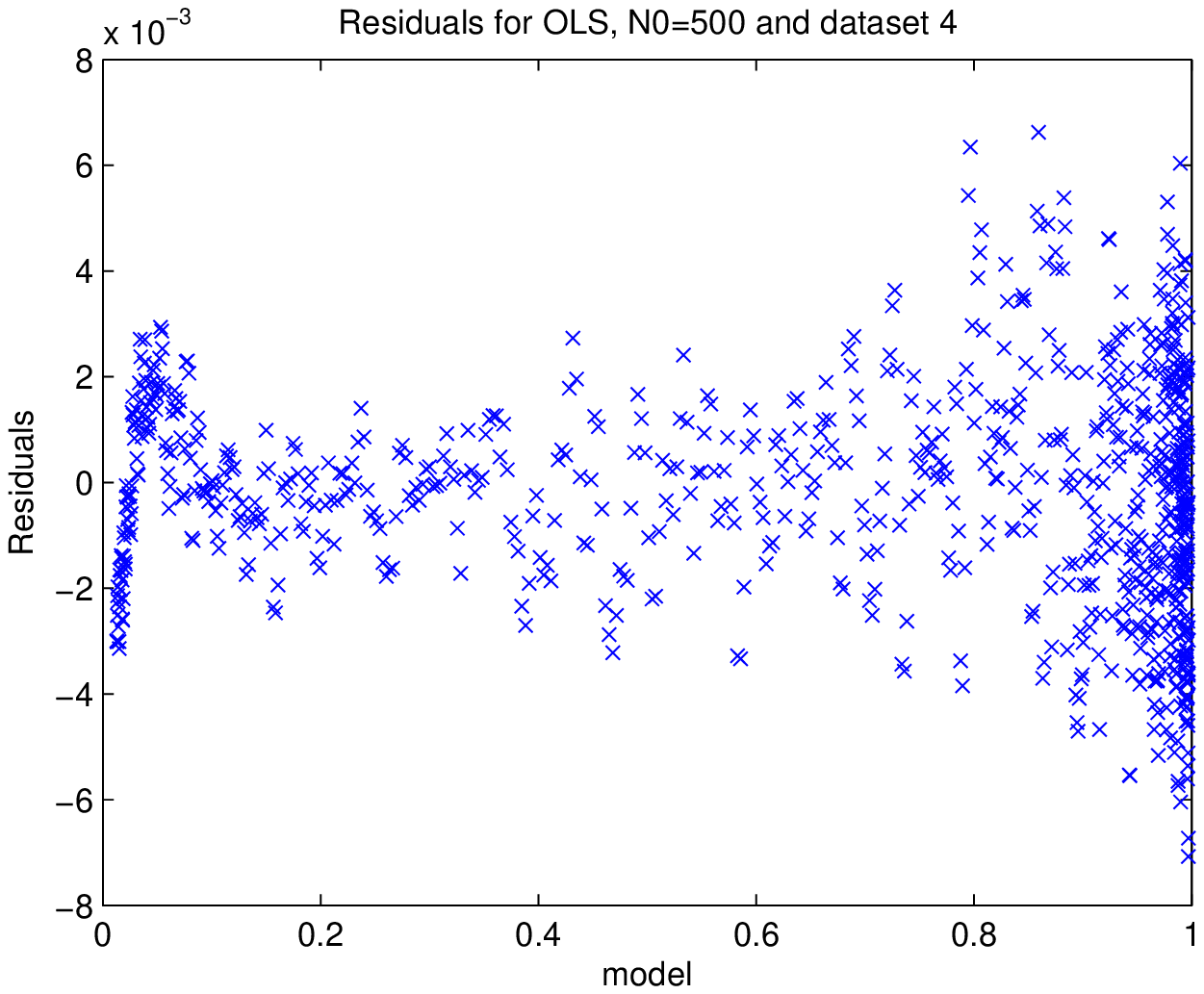}}\\
 (a) & (b)\\
 \end{tabular}
 \end{center}
 \caption{(a) $M(t_k)$ with OLS; (b) Residuals vs Model: OLS}
 \label{OLSdataset4}
 \end{figure}

 \begin{figure}[h!]
 \begin{center}
 \begin{tabular}{cc}
\resizebox{6.8cm}{!}{\includegraphics{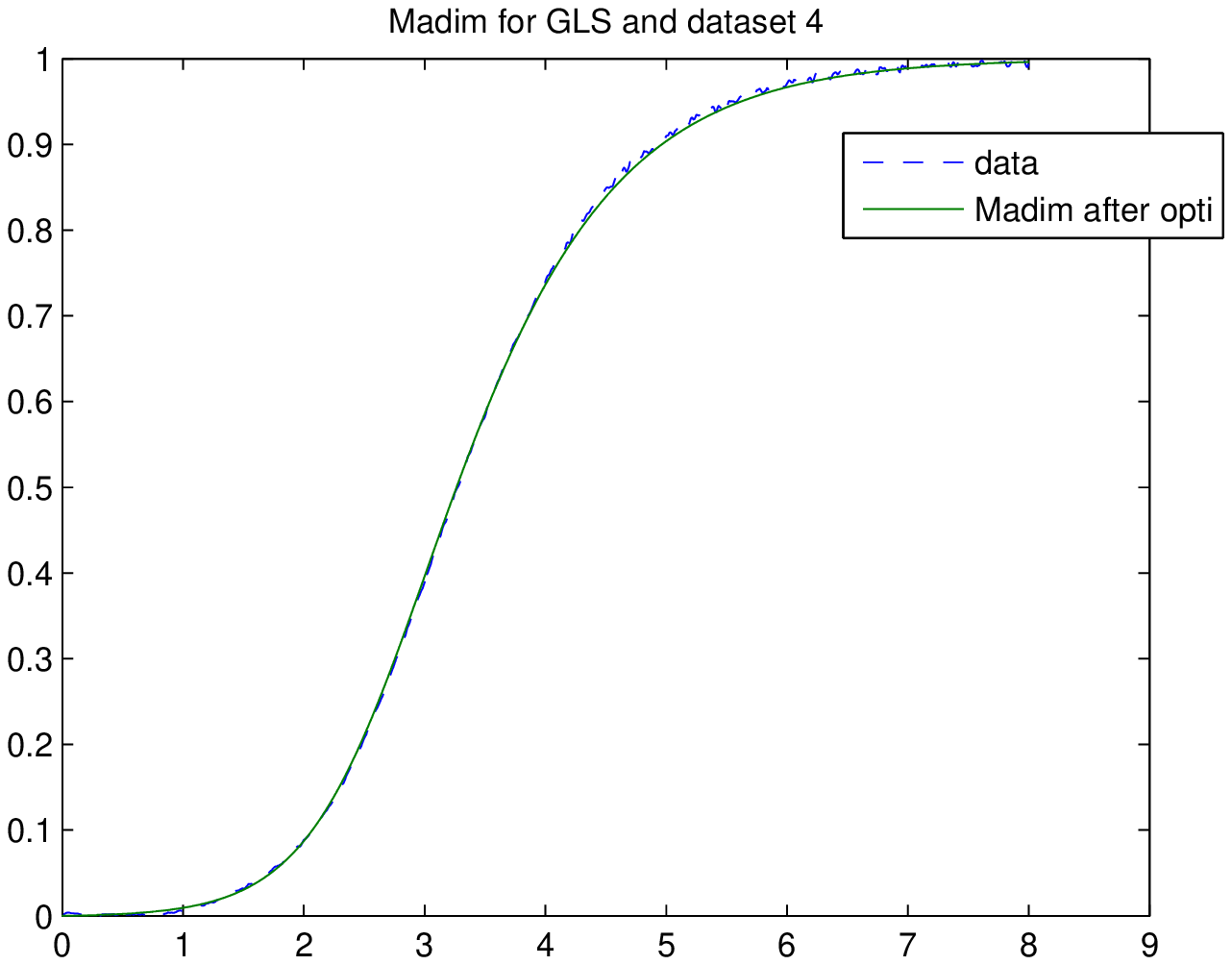}}&
\resizebox{6.8cm}{!}{\includegraphics{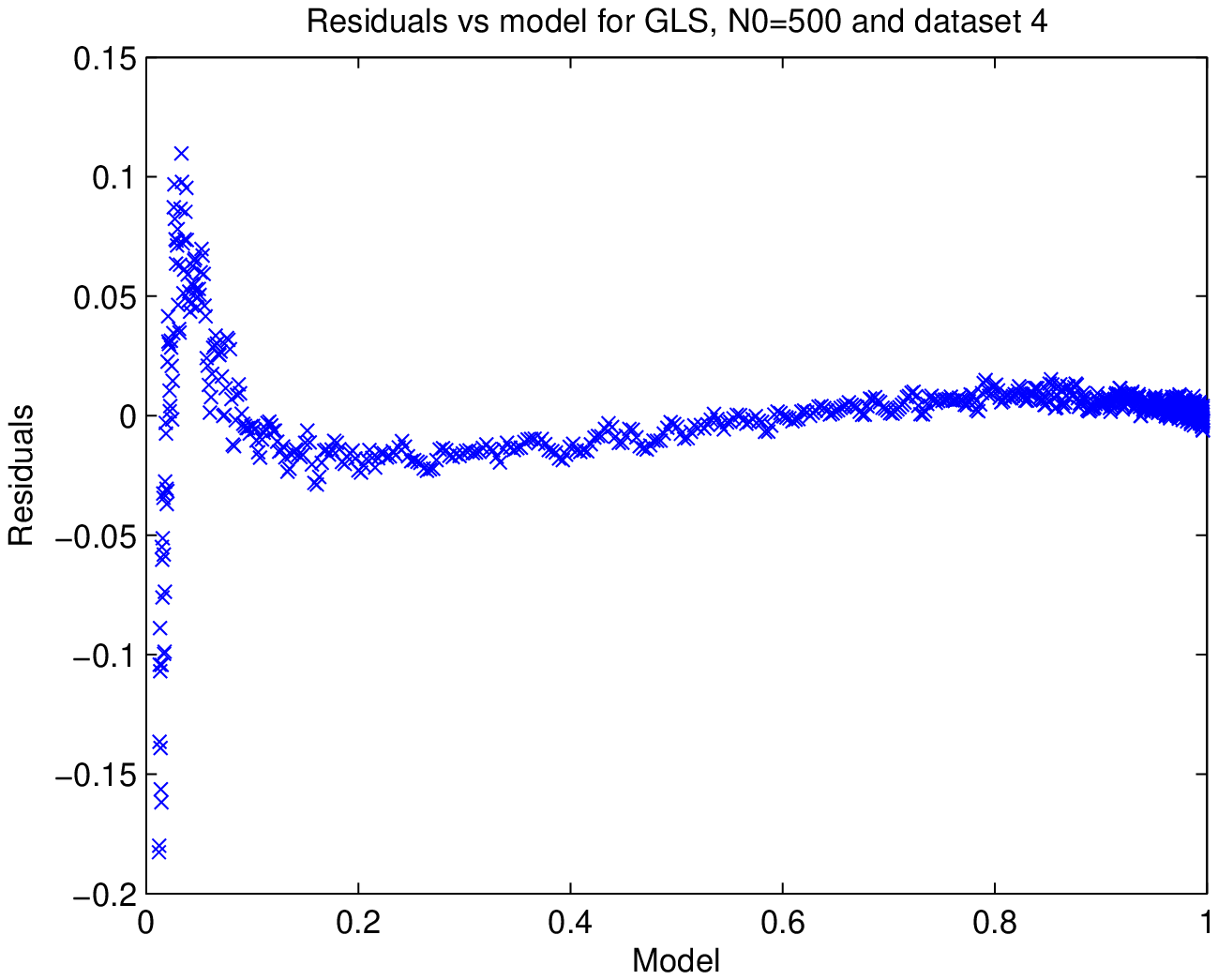}}\\
(a) & (b)\\
 \end{tabular}
 \end{center}
 \caption{ (a) $M(t_k)$ with GLS, $\gamma=1$; (b) Residuals vs Model: GLS}
 \label{GLSdataset4}
 \end{figure}

Based on these initial results and the speculation that early periods of the polymerization process may be somewhat stochastic in nature, we chose to subsequently use all the data sets on the intervals $[t_0,8]$ where $t_0$ is the first time when $M(t_0) >0.12$ (thus 12\% of the total polymerized mass). Moreover, we decided to use other values of $\gamma$ between 0 and 1 to test data set 4.

We thus carried out further investigations with inverse problems for data points $M(t_k) \geq 0.12$ and $i_0=2$ where we focused on the question of the most appropriate values of $\gamma$ to use in a generalized least squares approach (again see \cite{BHT2014} for further motivation and details). We then obtained the results with data set 4 depicted in Figure \ref{varying}.
 \begin{figure}[h!]
 \begin{center}
 \resizebox{16.cm}{!}{\includegraphics{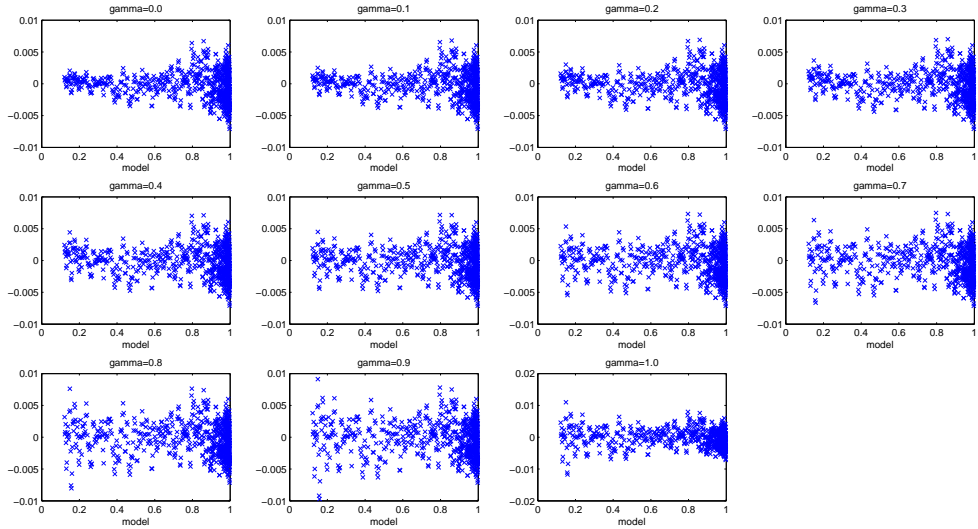}}
 \end{center}
 \caption{Residuals for data set 4 using different values of $\gamma$.}
 \label{varying}
 \end{figure}

 Analysis of these residuals suggest that either $\gamma=0.6$ or $\gamma=0.7$ might be satisfactory for use in a generalized least squares setting.


Motivated by these results, we next investigated the inverse problems for each of the four experimental data sets with initial concentration $c_0=200\mu$mol and $i_0=2$. We carried out the optimization over all data points with $M(t_k) \geq 0.12$ and used the  generalized least squares method with $\gamma=0.6$. The resulting graphics depicted in Figure \ref{4datasets} again suggest that $\gamma=0.6$ is a reasonable value to use in our subsequent analysis of the polyglutamine data with regard to its information content for inverse problem estimation and parameter uncertainty quantification.

 \begin{figure}[h]
 \begin{center}
 \resizebox{14cm}{!}{\includegraphics{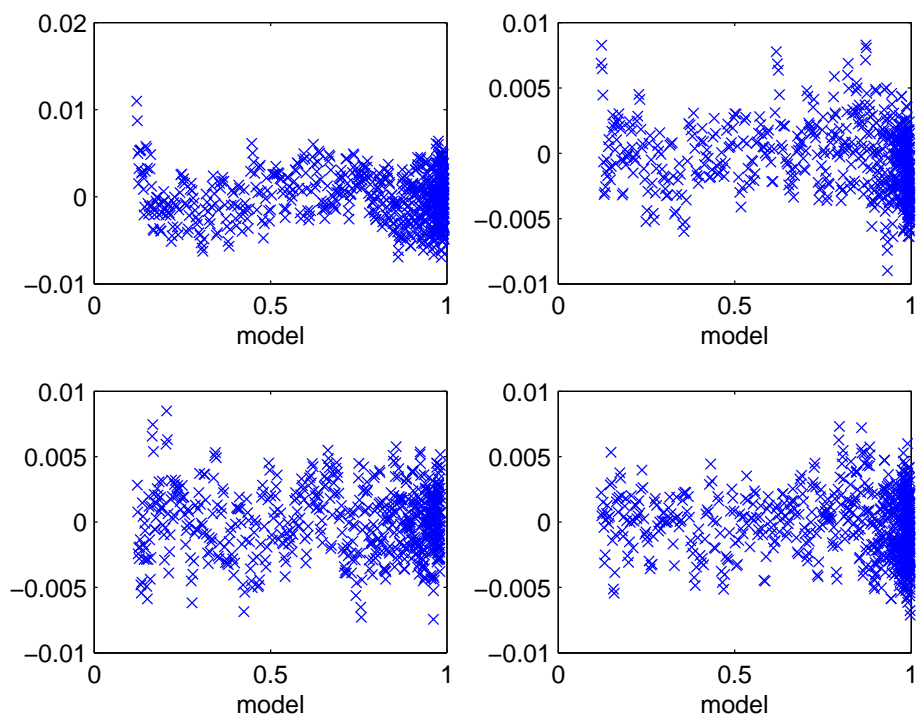}}
 \end{center}
 \caption{Residuals for the 4 experimental data sets using $\gamma=0.6$.}
 \label{4datasets}
 \end{figure}
\clearpage
\section{Standard Errors and Asymptotic Analysis}
\label{sec:SE}
\subsection{Standard Errors for Parameters Using GLS}

We employed first the asymptotic theory for parameter uncertainty summarized in \cite{BHT2014,BT2009,DG} and the references therein.
In the case of generalized least squares, the associated standard errors for the estimated parameters $\hat{\theta}=(k_I^+, ..., i_{max})$ (vector length $\kappa_{\theta}=9$) are given by the following construction (for details see Chap. 3.2.5 and 3.2.6 of \cite{BHT2014}):

Define the covariance matrix by the formula
\begin{align*}
SE_k=\sqrt{\Sigma_{kk}(\hat{\theta})},\hspace{0.5cm} k=1,...,9,
\end{align*}
where
\begin{align*}
\Sigma(\hat{\theta})=\hat{\sigma}^2 (\chi^T(\hat{\theta})W(\hat{\theta})\chi(\hat{\theta}))^{-1}.
\end{align*}
 Here $\chi$ is the sensitivity matrix of size $n\times \kappa_{\theta}$ ($n$ being the number of data points and $\kappa_{\theta}$ being the number of estimated parameters) and $W$ is defined by
\begin{align*}
W^{-1}(\hat{\theta})= \text{diag}(M(t_1;\hat{\theta})^{2\gamma},\dots,M(t_n;\hat{\theta})^{2\gamma}).
\end{align*}
We use the approximation of the variance
\begin{align*}
\sigma^2\approx\hat{\sigma}(\hat{\theta})^2=\frac{1}{n-\kappa_{\theta}}\sum_{i=1}^n \frac{1}{M(t_i;\hat{\theta})^{2\gamma}}(M(t_i,\hat{\theta})- y_i)^2.
\end{align*}

To obtain a finite standard error using asymptotic theory, the $9\times 9$ matrix  $F=\chi^T(\hat{\theta})W(\hat{\theta})\chi(\hat{\theta})$ thus must be invertible.
In the above problem we do indeed obtain a good fit of the curve and good residuals (for the sake of brevity, not depicted here!). However, we also found that the condition number of the matrix $$F=\chi^T(\hat{\theta})W(\hat{\theta})\chi(\hat{\theta})$$ is $\kappa=10^{24}$. Looking more closely at the matrix $F$ reveals a near linear dependence between certain rows, hence the large condition number. We thus quickly reach the {\em following conclusions}:

\begin{enumerate}
	\item We obtain a set of parameters for which the model fits well, but we {\em cannot} have any reasonable confidence in them using the asymptotic theories from statistics e.g., see the references given above.
	\item We suspect that it may  not be possible to obtain sufficient information from our data set curves to estimate all 9 parameters with a high degree of confidence! This is based on our calculations with the corresponding Fisher matrices as well our prior knowledge in that the graphs depicted in Figure \ref{datasets} are very similar to Logistic or Gompertz curves which can be quite well fit with parameterized  models with only 2 or 3 carefully chosen parameters!
\end{enumerate}
To assist in initial understanding of these issues, we consider the associated sensitivity matrices $\chi=\frac{\partial M}{\partial \theta}$.

\subsection{Sensitivity Analysis}
For the sensitivity analysis, we follow \cite{BHT2014,BT2009}. Hereafter all our analysis will be carried using data set 4 and the best estimate  $\hat{\theta}$ obtained for the latter. We find that the model is sensitive mainly to four parameters: $k_I^+, k_I^-, k_{on}^{N}, k_{off}^{N}$. The sensitivities for the remaining parameters are on an order of magnitude of  $10^{-6}$ or less.
It also shows some sensitivity with respect to $x_1$. However, the parameter $x_1$ appears in the model only as factor $x_1 i_{max}$. The sensitivities depicted below use $\hat{\theta}$ for the nine best fit GLS parameters , i.e., $\hat{\theta}$ for  $\kappa_{\theta}=9$.

\begin{figure}[h!]
 \begin{tabular}{cc}
 \resizebox{7.0cm}{!}{\includegraphics{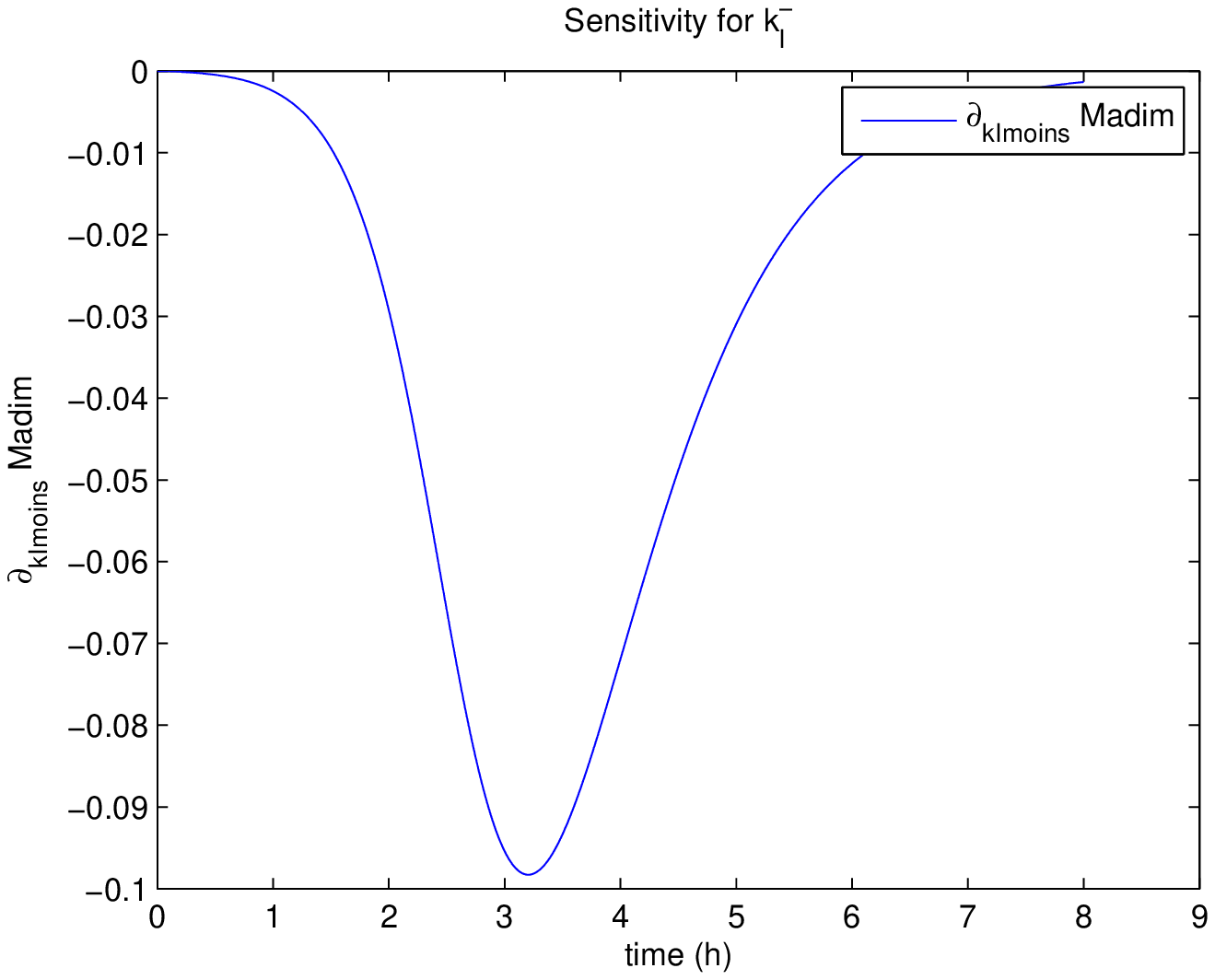}}&
 \resizebox{7.0cm}{!}{\includegraphics{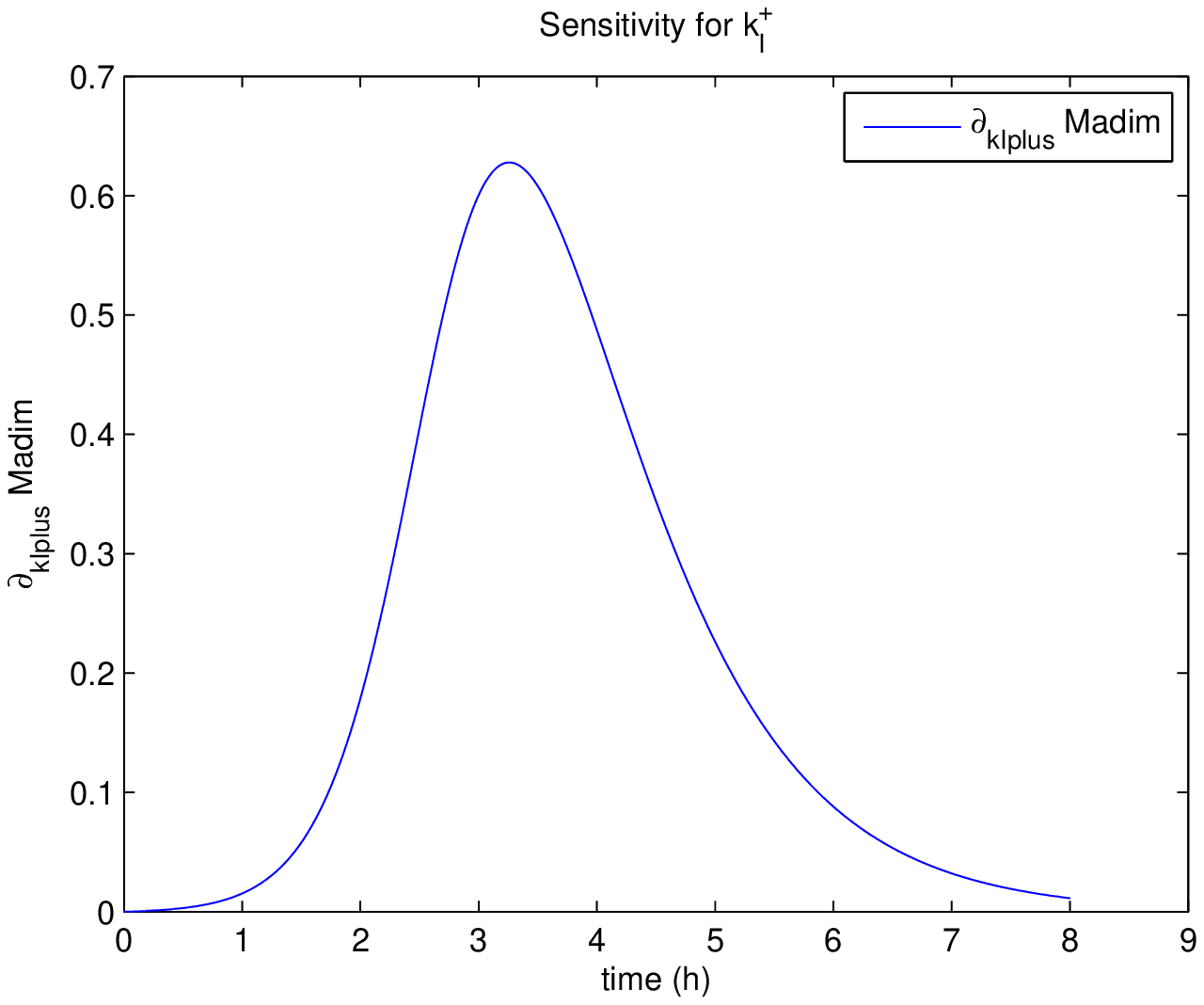}}\\
 (a) & (b)\\
 \end{tabular}
\caption{(a) Sensitivity w.r.t. $k_I^{-}$; (b) Sensitivity w.r.t. $k_I^{+}$}
 \end{figure}

\begin{figure}[h!]
 \begin{tabular}{cc}
 \resizebox{7.0cm}{!}{\includegraphics{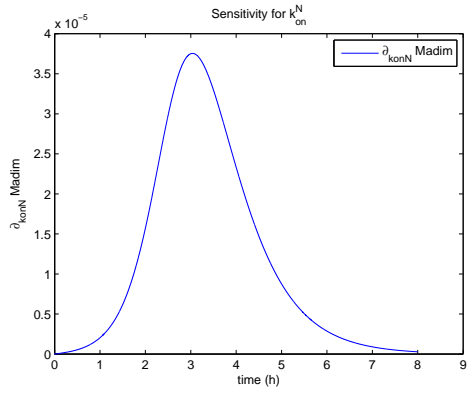}}&
 \resizebox{7.0cm}{!}{\includegraphics{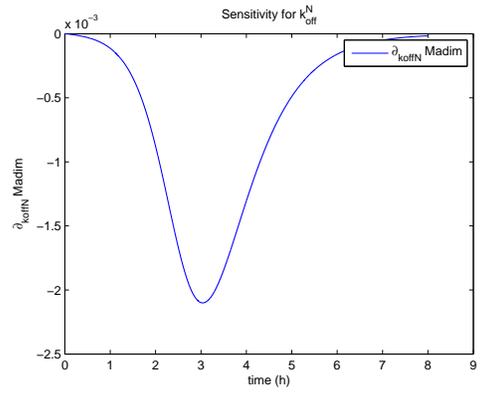}} \\
 (a) & (b)\\
 \end{tabular}
\caption{(a) Sensitivity w.r.t. $k_{on}^{N}$; (b) Sensitivity w.r.t. $k_{off}^{N}$ }
 \end{figure}


\begin{figure}[h!]
 \begin{tabular}{cc}
 \resizebox{7.0cm}{!}{\includegraphics{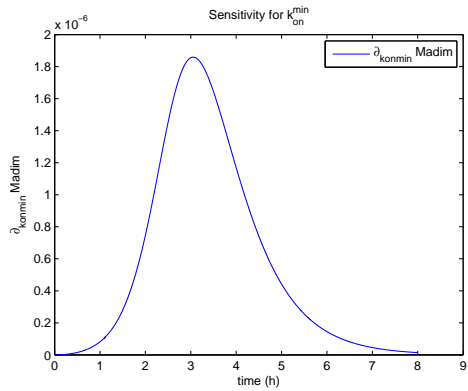}}&
 \resizebox{7.0cm}{!}{\includegraphics{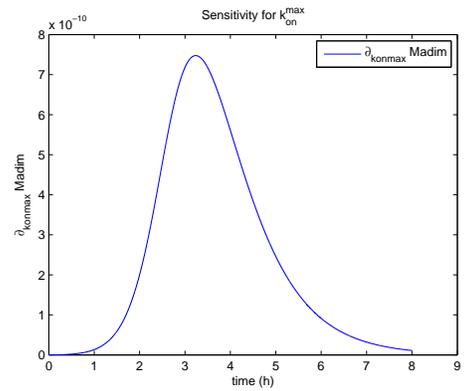}} \\
 (a) & (b)\\
 \end{tabular}
\caption{(a) Sensitivity w.r.t. $k_{on}^{min}$; (b) Sensitivity w.r.t. $k_{off}^{max}$ }
 \end{figure}

\begin{figure}[h!]
 \begin{tabular}{cc}
 \resizebox{7.0cm}{!}{\includegraphics{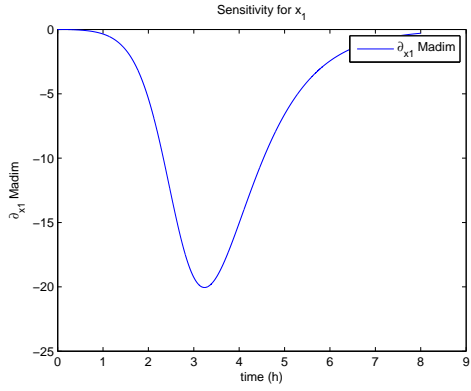}}&
 \resizebox{7.0cm}{!}{\includegraphics{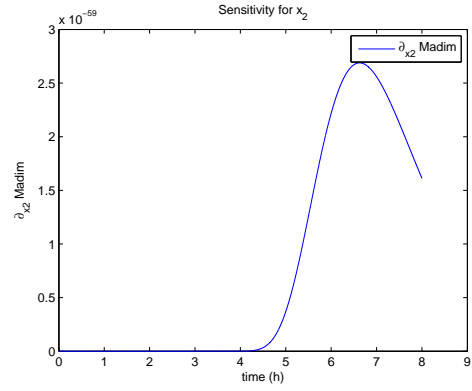}} \\
 (a) & (b)\\
 \end{tabular}
\caption{(a) Sensitivity w.r.t. $x_1$; (b) Sensitivity w.r.t. $x_2$ }
 \end{figure}
\begin{figure}[h!]
 \begin{tabular}{cc}
 \resizebox{7.0cm}{!}{\includegraphics{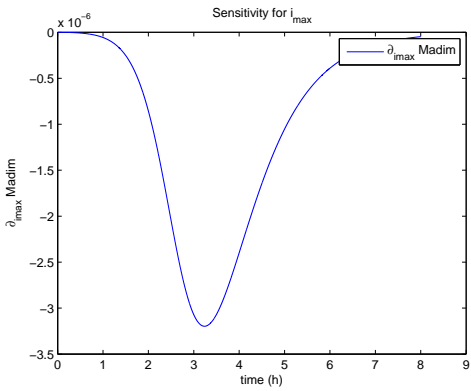}}&
 \resizebox{7.0cm}{!}{\includegraphics{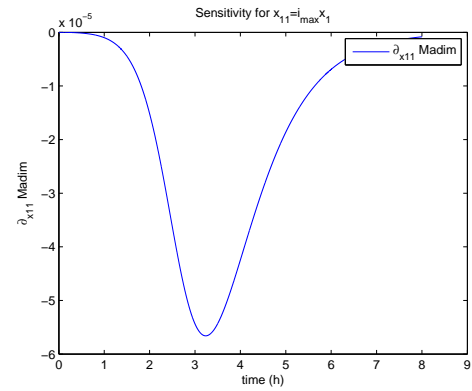}} \\
 (a) & (b)\\
 \end{tabular}
\caption{(a) Sensitivity w.r.t. $i_{max}$; (b) Sensitivity w.r.t. $x_{11}=i_{max}x_1$ }
 \end{figure}


\clearpage

\section{Sensitivity Motivated Inverse Problems}\label{sec:boot}
Based on the sensitivity findings depicted above, we investigated a series of inverse problems in which we attempted to estimate an increasing number of parameters
beginning first with the fundamental parameters $k_I^+$ and $k_I^-$. In each of these inverse problems we attempted to ascertain uncertainty bounds  for the estimated parameters using both the asymptotic theory described above and a generalized least squares version of bootstrapping \cite{CR,CWR,MDnotes,DiCiEfron,Efronbook}.

A quick outline of the appropriate bootstrapping algorithm is given next.


\subsection{Bootstrapping Algorithm: Nonconstant Variance Data}
We suppose now that we are given experimental data
$(t_1,y_1),\ldots,(t_n,y_n)$
 from the underlying observation process
\begin{equation}\label{relerr}
Y_i=M(t_i;\theta_0)+  M(t_i;\theta_0)^{\gamma}\widetilde{\mathcal{E}_i},
\end{equation}
where $i=1,\dots,n$ and the $\widetilde{\mathcal{E}_i}$ are \textit{i.i.d.} with mean zero and constant variance $\sigma_0^2$.  Then we see
that $\mathbb{E}(Y_i) =M(t_i;\theta_0)$ and ${Var}(Y_i) =\sigma_0^2
M^{2\gamma}(t_i,\theta_0)$, with associated corresponding realizations of
$Y_i$ given by
\begin{equation*}
y_i=M(t_i;\theta_0)+M(t_i;\theta_0)^{\gamma}\widetilde{\epsilon_i}.
\end{equation*}

A standard algorithm  can
be used to compute the corresponding {\em bootstrapping estimate}
$\hat{\theta}_{boot}$ of $\theta_0$ and its empirical
distribution. We treat the general case for nonlinear dependence of the model output on the parameters $\theta$.
The algorithm is given as follows.

\begin{enumerate}
\item First obtain the estimate
$\hat{\theta}^0$ from the entire
sample $\{y_i\}$ using the GLS given in \eqref{WLSform} with $\gamma=1$. An estimate $\hat{\theta}_{boot}$ can be solved for iteratively as follows.

\item Define the nonconstant variance standardized residuals
\begin{equation*}
 \bar{s}_i = \frac{y_i -
 M(t_;\hat{\theta}^0)}{M(t_i;\hat{\theta}^0)^{\gamma}}, \quad i=1,2,\dots,n.
\end{equation*}

Set $m=0$.
\item Create a bootstrapping sample of size $n$ using random sampling with
replacement from the data (realizations)
\{$\bar{s}_1$,\ldots,$\bar{s}_n$\} to form a bootstrapping sample
$\{s_1^m,\ldots,s_n^m\}$.

\item Create bootstrapping sample points
\begin{equation*}
 y_i^m =M(t_i;\hat{\theta}^0)+M(t_i;\hat{\theta}^0)^{\gamma}s_i^m,
\end{equation*}
where $i=1$,\ldots,$n$.

\item Obtain a new estimate
$\hat{\theta}^{m+1}$
from the bootstrapping sample $\{y_i^m\}$ using GLS.

\item Set $m=m+1$
and repeat steps 3--5 until $m\geq M$  where M is large (e.g., M=1000). 

\end{enumerate}

We then calculate the mean, standard
error, and confidence intervals using the
formulae
\begin{eqnarray}\label{bootform}
&\hat{\theta}_{boot} = \frac{1}{M}\sum_{m=1}^M \hat{\theta}^m, \nonumber\\
&{\text Var}({\theta}_{boot}) =\frac{1}{M-1} \sum_{m=1}^M (\hat{\theta}^m-\hat{\theta}_{boot})^T(\hat{\theta}^m-\hat{\theta}_{boot}),\\
& \textrm{SE}_k(\hat{\theta}_{boot}) =
\sqrt{Var({\theta}_{boot})_{kk}}. \nonumber
\end{eqnarray}
where ${\theta}_{boot}$ denotes the bootstrapping estimator.

\subsection{Estimation of two parameters}

We first carried out estimation for the 2 parameters $k_{I}^+$ and $k_I^{-}$. We use the GLS formulation with $\gamma=0.6$. We fix globally (based on previous estimations with DS 4) the parameter values
\begin{align*}
\begin{array}{|c|c |c |c |c|c|c|}
\hline
k_{on}^N & k_{off}^N & k_{on}^{min}& k_{on}^{max}& x_1 & x_2 & i_{max} \\
\hline
4616.962 &93.332 & 1684.381 & 1.5152 \cdot 10^9 &  0.0626 & 0.859 & 3.542\cdot 10^5\\
\hline
\end{array}
\end{align*}
and used the initial guesses for the parameters given by
\begin{align*}
\begin{array}{|c |c |c |}
\hline
 &k_I^+ & k_I^{-} \\
\hline
q_0 & 2.1600 & 10.9270 \\
\hline
\end{array}
\end{align*}
%

We then used the bootstrapping algorithm and obtained the following means and standard errors for $M=1000$ which, as reported below, compare quite well with the asymptotic theory estimates. The corresponding distributions are shown in Figures  \ref{fig:boot1000_2p_kip} and \ref{fig:boot1000_2p_kim}.

\begin{align*}
\begin{array}{|c |c |c |c|c|c|c|c|}
\hline
 &k_I^+ (boot) (GLS) & k_I^{-} (boot) (GLS)&  k_I^+ (asymp) (GLS)& k_I^- (asymp) (GLS)\\
\hline
mean & 2.158 & 10.911  & 2.157  &  10.911 \\
\hline
SE &  0.0044  &  0.0247   &  0.00396 &  0.0225 \\
\hline
\end{array}
\end{align*}

\begin{figure}[htp]
  \centering
	\resizebox{8cm}{!}{\includegraphics{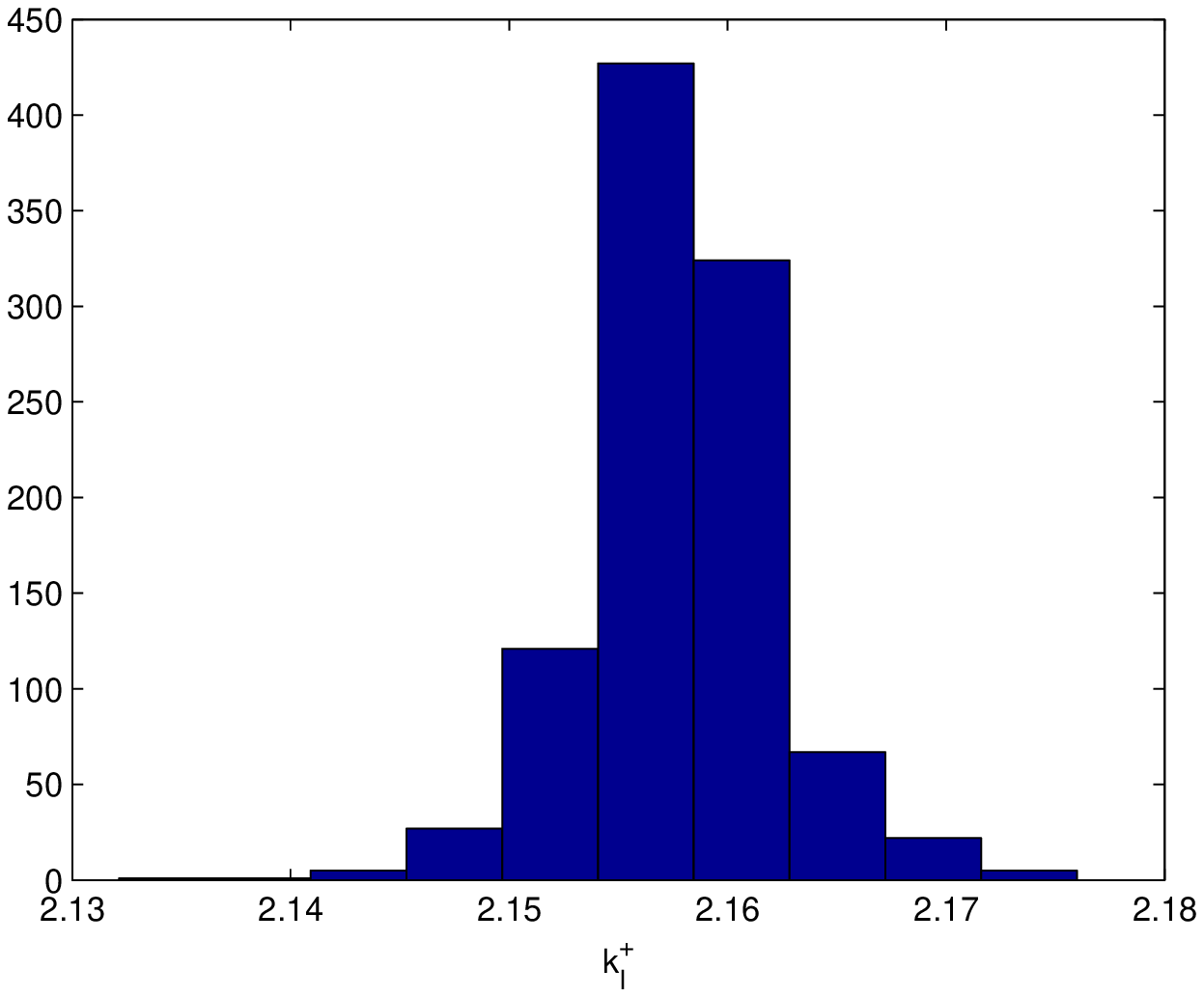}}
	\caption{Two parameters estimation ($k_I^+$, $k_I^-$). Bootstrapping distribution for $k_I^+$. We use GLS and M=1000 runs.}	
\label{fig:boot1000_2p_kip}
\end{figure}

\begin{figure}[htp]
  \centering
	\resizebox{8cm}{!}{\includegraphics{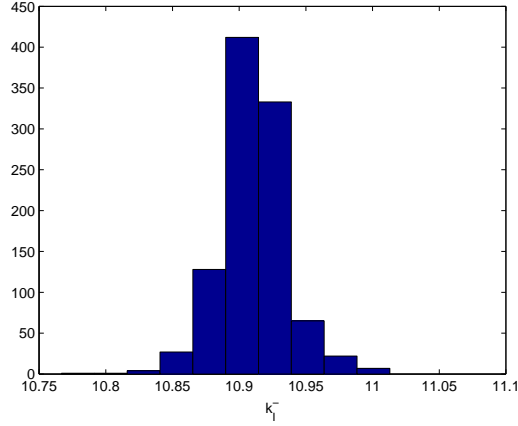}}
	  \caption{Two parameters estimation ($k_I^+$, $k_I^-$). Bootstrapping distribution for $k_I^-$. We use GLS and M=1000 runs.}	
	\label{fig:boot1000_2p_kim}
\end{figure}

\subsection{GLS Estimation of 3 Parameters}

 We tried next to estimate 3 parameters. We again used the GLS formulation with $\gamma=0.6$. Once again we fixed all the parameters describing the domain and the polymerization function $k_{on}$ and we also fixed either $k_{off}^N$ or $k_{on}^N$ in the corresponding inverse problems.

\subsection{GLS Estimation for $k^+_I, \;k^-_I$ and $k_{on}^N$}\label{sec:WLSestkonn}
We fixed values as follows:
\begin{align*}
\begin{array}{|c |c |c |c|c|c|}
\hline
k_{off}^N & k_{on}^{min}& k_{on}^{max}& x_1 & x_2 & i_{max} \\
\hline
93.33 & 1684.38 & 1.5 \cdot 10^9 &  0.062 & 0.859 & 3.5\cdot 10^5\\
\hline
\end{array}
\end{align*}
We used as initial parameter values:
\begin{align*}
\begin{array}{|c |c |c |c|}
\hline
 &k_I^+ & k_I^{-}& k_{on}^{N} \\
\hline
q_0 & 2.1600 & 10.9270 & 4616.962\\
\hline
\end{array}
\end{align*}
We obtained the estimated parameters together with the corresponding standard errors, variances and the condition numbers $\kappa$ of the corresponding sensitivity matrices for the four data sets as reported below. The $95\%$ confidence results based on the  asymptotic theory are also depicted for DS 4 in Figure~\ref{fig:err_bar_pl}.

\begin{align*}
\begin{array}{|c |c |c |c|c|c|c|}
\hline
 &k_I^+ & k_I^{-}& k_{on}^{N} & SE & \sigma^2 & \kappa \\
\hline
DS 1 & 2.26 & 13.49  & 4616.96& ( .012,   .099,  53.925) &  8.52\cdot 10^{-6} &  8.89 \cdot 10^{10} \\
\hline
DS 2 & 2.99 &  16.20 & 4616.96 &  (.021, .151,  56.691) &  9.67\cdot 10^{-6}  & 4.37 \cdot 10^{10} \\
\hline
DS 3 & 2.18 & 15.76  & 9840.31& ( .011,   .103,  90.466)&  6.45\cdot 10^{-6} &   3.94\cdot 10^{11}\\
\hline
DS 4 &2.16&10.91&4616.96& ( 0.0089,   0.0649,  45.262)  &  6.36\cdot 10^{-6}&   7.14 \cdot 10^{10}\\
\hline
\end{array}
\end{align*}

\begin{figure}[htp]
  \centering
	\resizebox{10cm}{!}{\includegraphics{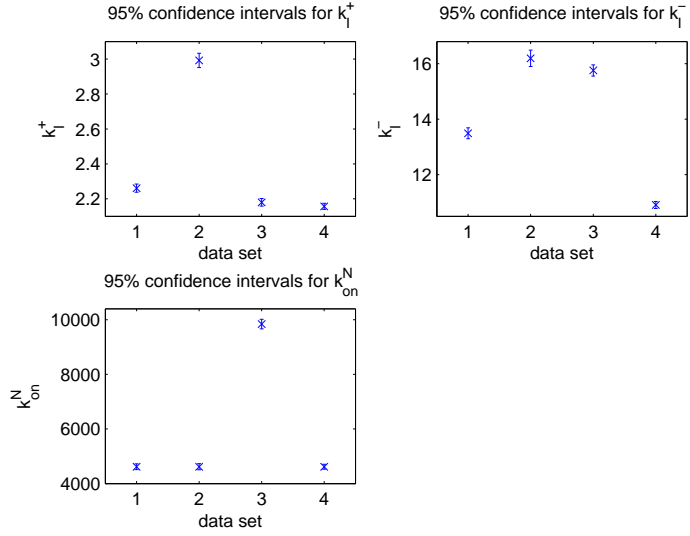}}
\caption{
Confidence Intervals}
\label{fig:err_bar_pl}
\end{figure}

To compare these asymptotic results with bootstrapping, we carried out bootstrapping with Data Set (DS) 4 for the estimation of $k_{I}^+$, $k_I^{-}$ and $k_{on}^N$ with the same initial values as above.
We then obtained the following means and standard errors for a run with 
$M=1000$, in comparison to the asymptotic theory.
\begin{align*}
\begin{array}{|c |c |c |c|c|c|c|c|c|c|}
\hline
 &k_I^+ (boot) & k_I^{-} (boot)& k_{on}^{N} (boot) & k_I^+ (asymp) & k_I^- (asymp) &  k_{on}^N (asymp)\\
\hline
mean & 2.153 &  10.887 & 4616.962  & 2.157 & 10.910 & 4616.962  \\
\hline
SE &  0.0039 &  0.0219  &   0.00003 & 0.0089 & 0.0649& 45.262 \\
\hline
\end{array}
\end{align*}

Of particular interest are the values obtained for $k^N_{on}$ and the bootstrapping standard errors for $k_{on}^N$ which are extremely small. It should be noted that the sensitivity of the model output on $k_{on}^N$ is also very small. Thus one might conjecture that the iterations in the bootstrapping algorithm do not change the values of $k^N_{on}$ very much and hence one observes the extremely small SE that are produced for the bootstrapping estimates.

\begin{figure}[htp]
  \centering
	\resizebox{8cm}{!}{\includegraphics{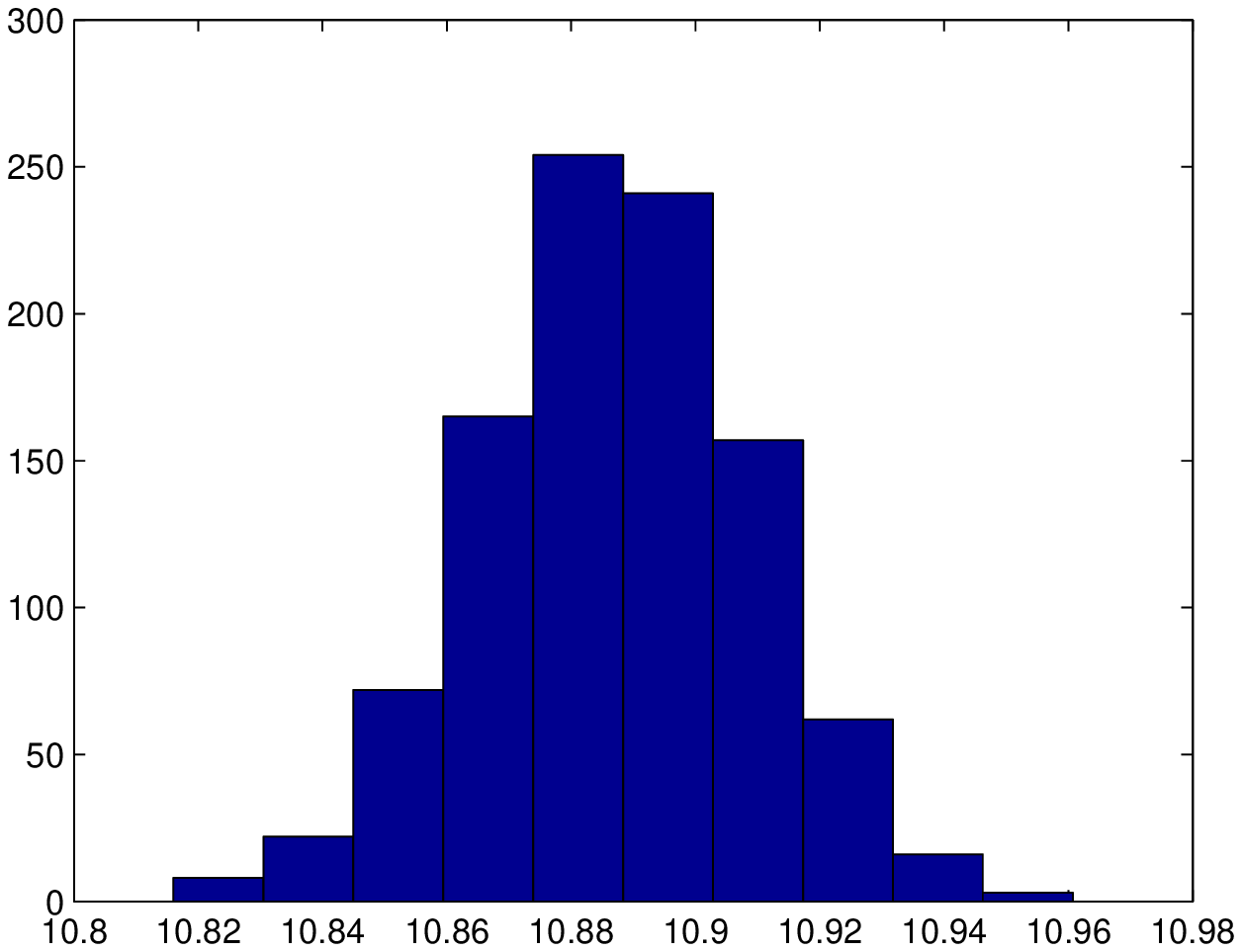}}
	  \caption{Estimation for $k_{I}^+$, $k_I^{-}$ and $k_{on}^N$: Bootstrapping distribution  for $k_I^-$ for GLS and 1000 runs.}	
	\label{fig:boot400_OLS_kim}
\end{figure}
\begin{figure}[htp]
  \centering
	\resizebox{8cm}{!}{\includegraphics{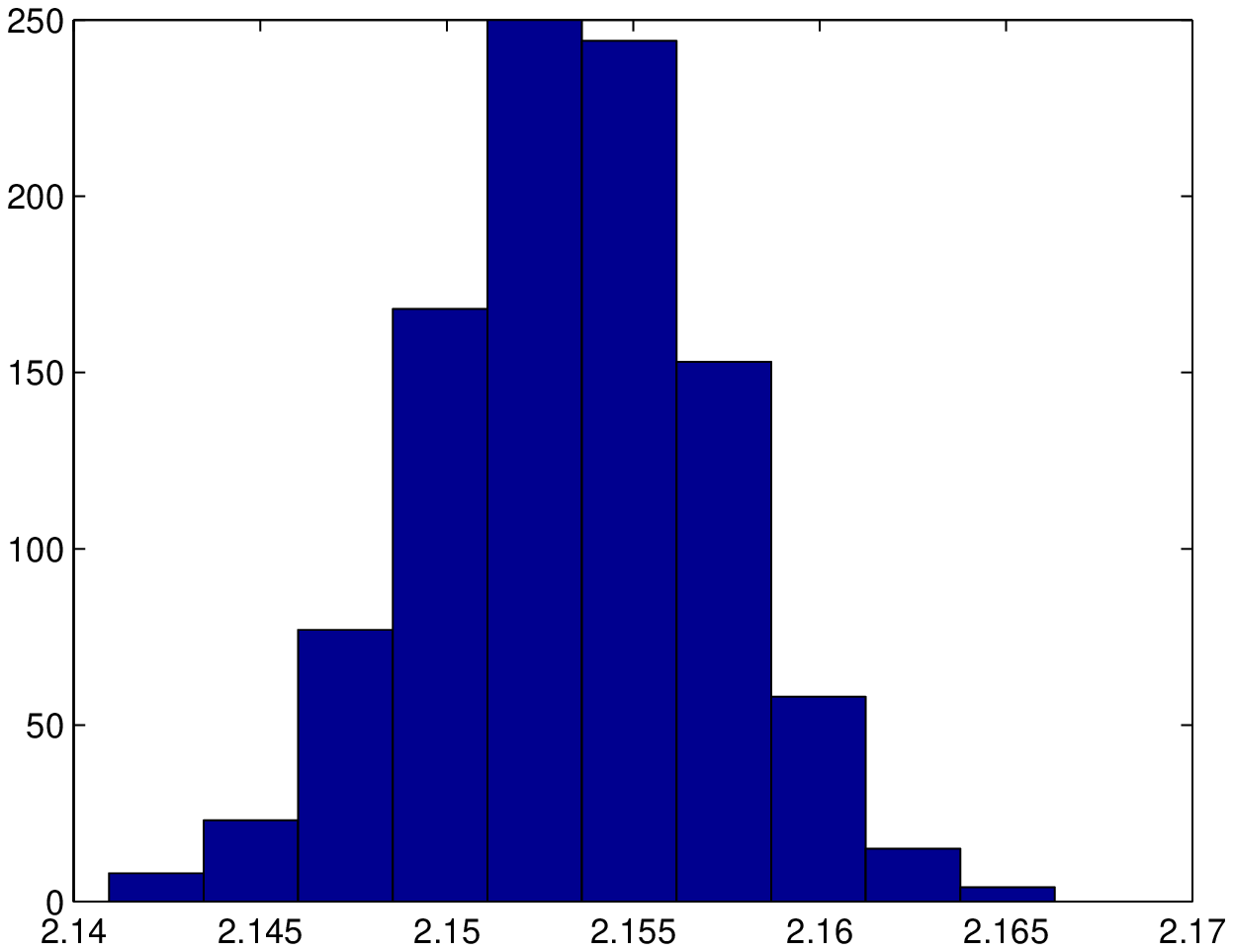}}
	\caption{Estimation for $k_{I}^+$, $k_I^{-}$ and $k_{on}^N$: Bootstrapping distribution for $k_I^+$ for GLS and 1000 runs.}	
\label{fig:boot400_OLS_kip}
\end{figure}
\begin{figure}[htp]
  \centering
	\resizebox{8cm}{!}{\includegraphics{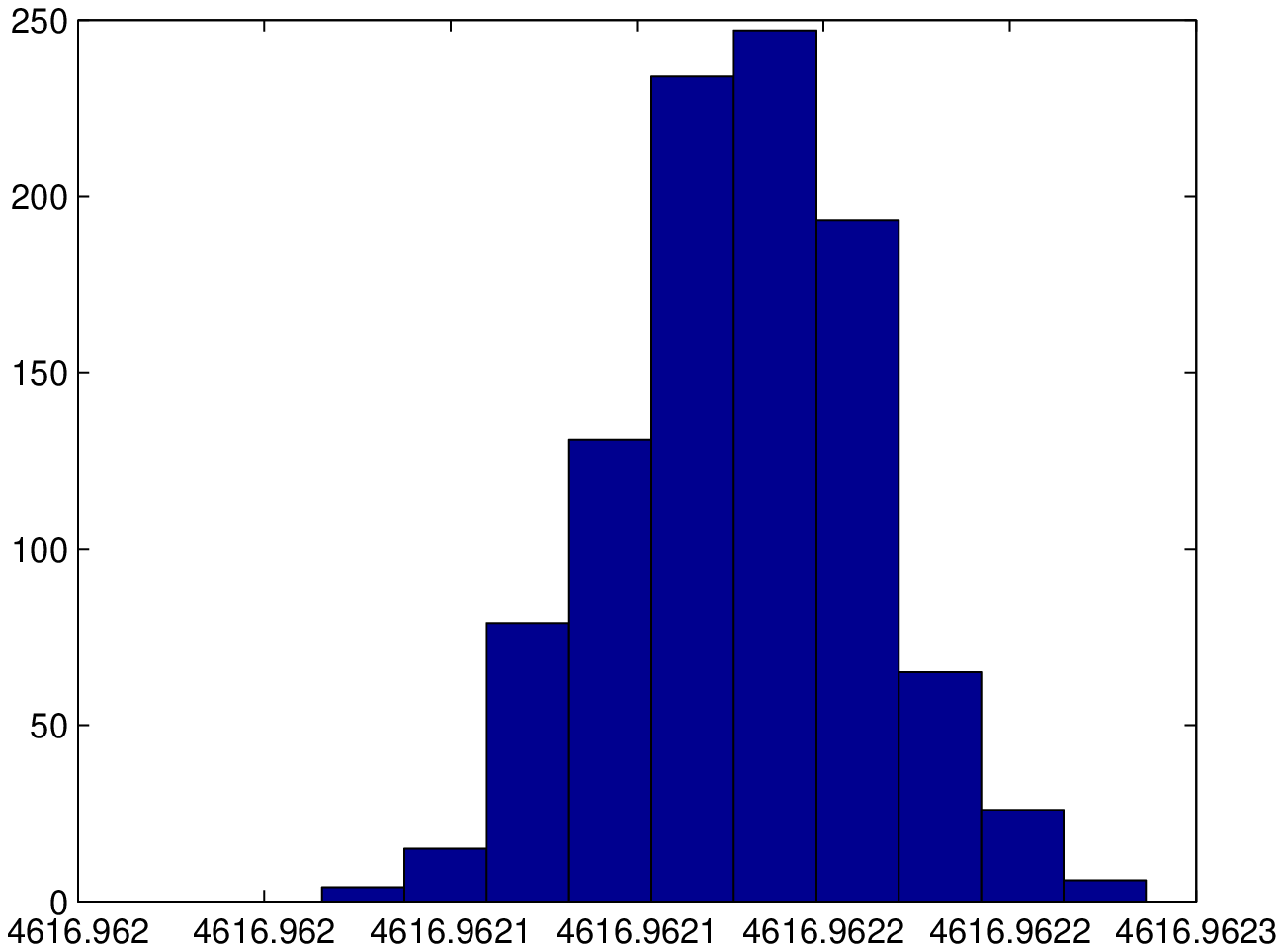}}
	\caption{Estimation for $k_{I}^+$, $k_I^{-}$ and $k_{on}^N$: Bootstrapping distribution for $k_{on}^N$ for GLS and 1000 runs.}	
\label{fig:boot400_OLS_kip}
\end{figure}

\clearpage
\subsection{GLS estimation for $k^+_I, K^-_I$ and  $k_{off}^N$ }

In another test, we fixed $k_{on}^N$ and instead estimate $k_{off}^N$ (along with $k^+_I$ and $k^-_I$). We use the fixed values:
\begin{align*}
\begin{array}{|c |c |c |c|c|c|}
\hline
k_{on}^N & k_{on}^{min}& k_{on}^{max}& x_1 & x_2 & i_{max} \\
\hline
4616.962 & 1684.381 & 1.5152 \cdot 10^9 &  0.0626 & 0.859 & 3.542\cdot 10^5\\
\hline
\end{array}
\end{align*}
and the initial guesses for the parameters to be estimated given by:

\begin{align*}
\begin{array}{|c |c |c |c|}
\hline
 &k_I^+ & k_I^{-}& k_{off}^{N} \\
\hline
q_0 & 2.1600 & 10.9270 & 108.256\\
\hline
\end{array}
\end{align*}
We obtained the estimated parameters and corresponding SE.

\begin{align*}
\begin{array}{|c |c |c |c|c|c|c|}
\hline
 &k_I^+ & k_I^{-}& k_{off}^{N} & SE & \sigma^2 & \kappa \\
\hline
DS 1 &  2.203 &  12.997  & 99.861& (  .011, .091, 1.208) &  8.165\cdot 10^{-6} &  4.912 \cdot 10^{7} \\
\hline
DS 2 & 2.893 & 15.474  &  100.019 &  ( .019, .137, 1.279) &   9.323 \cdot 10^{-6}  &    2.486 \cdot 10^{7} \\
\hline
DS 3 & 2.168 & 15.631  & 41.935 & ( .011, .102, 0.424)&  6.435 \cdot 10^{-6} & 9.125 \cdot 10^{6}\\
\hline
DS 4 &2.181& 11.090& 90.536 & (.009, .066 , 0.936)  &  6.289 \cdot 10^{-6}&   3.043 \cdot 10^{7}\\
\hline
\end{array}
\end{align*}


Also in this case, we carried out  bootstrapping for DS 4. The bootstrapping distributions for $k_I^+$, $k_I^-$ and $k_{off}^N$ are found in Figures \ref{fig:boot_1000_3p_koffn_kip}-\ref{fig:boot_1000_3p_koffn_koffn}. We then obtained the following means and standard errors for a run with $M=1000$ in comparison to the asymptotic theory.
\begin{align*}
\begin{array}{|c |c |c |c|c|c|c|c|c|c|}
\hline
 &k_I^+ (boot) & k_I^{-} (boot)& k_{off}^{N} (boot) & k_I^+ (asymp) & k_I^- (asymp) &  k_{off}^N (asymp)\\
\hline
mean  &2.169 & 11.013&  91.254 &2.181& 11.090& 90.536  \\
\hline
SE &  0.0094 &  0.0699 &  1.0392  & 0.009 & 0.066 & 0.936\\
\hline
\end{array}
\end{align*}

\begin{figure}[htp]
  \centering
	\resizebox{8cm}{!}{\includegraphics{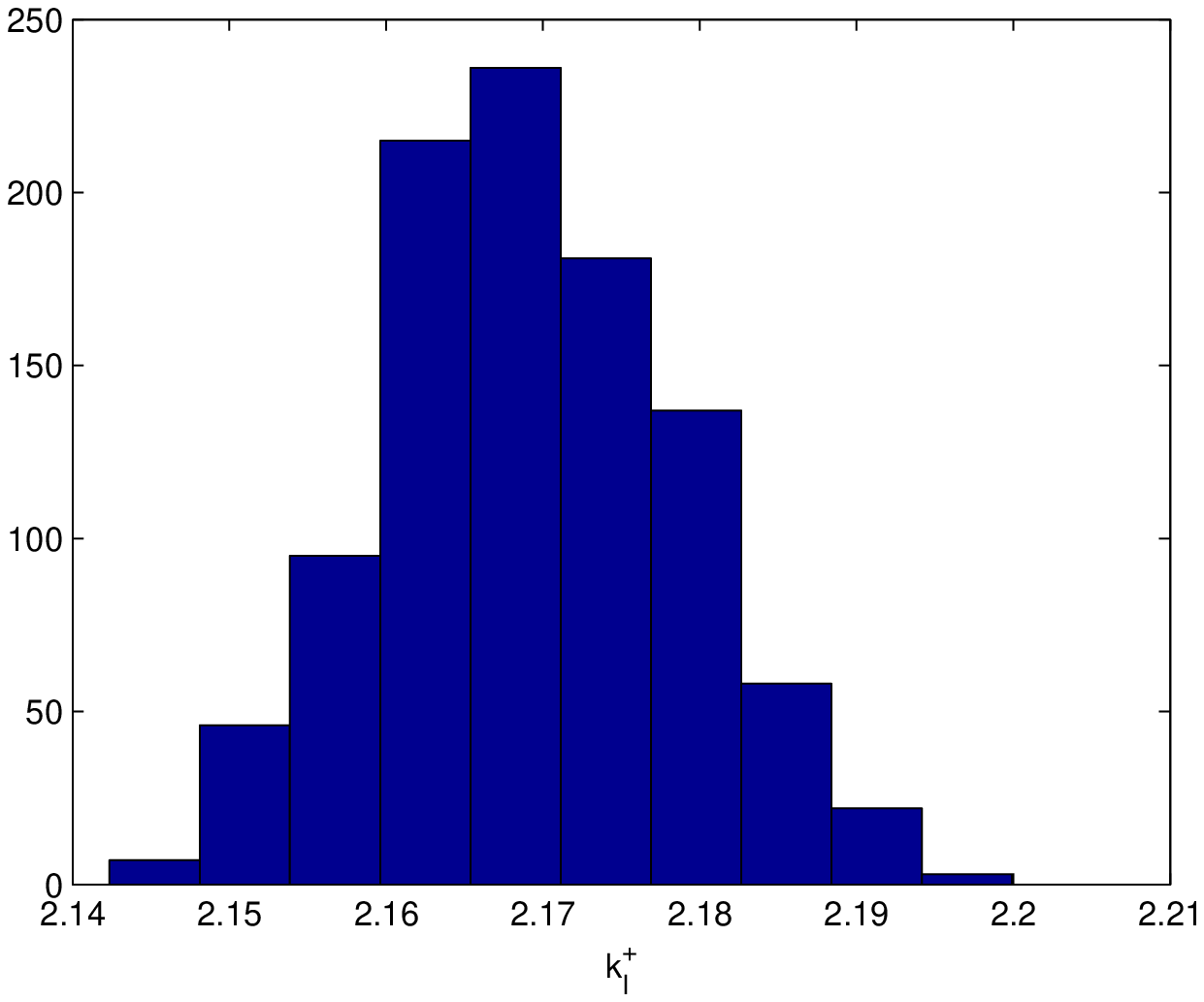}}
	\caption{Three parameters estimation ($k_I^+$, $k_I^-$ and $k_{off}^N$): Bootstrapping distribution for $k_I^+$. We used GLS and M=1000 runs.}	
\label{fig:boot_1000_3p_koffn_kip}
\end{figure}
\begin{figure}[htp]
  \centering
	\resizebox{8cm}{!}{\includegraphics{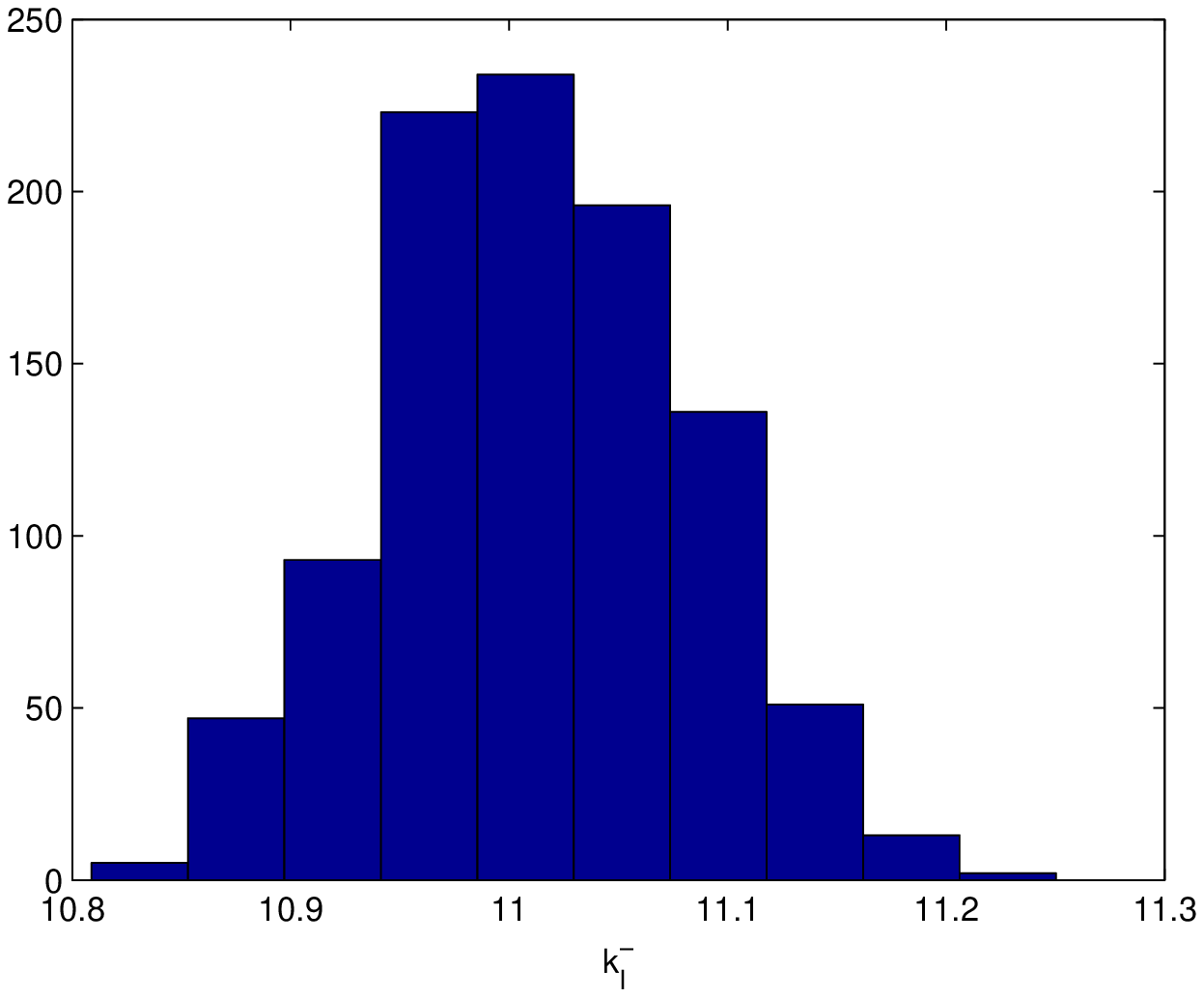}}
	  \caption{Three parameters estimation ($k_I^+$, $k_I^-$ and $k_{off}^N$): Bootstrapping distribution for $k_I^-$. We used GLS and M=1000 runs.}
	\label{fig:boot_1000_3p_koffn_kim}
\end{figure}

\begin{figure}[htp]
  \centering
	\resizebox{8cm}{!}{\includegraphics{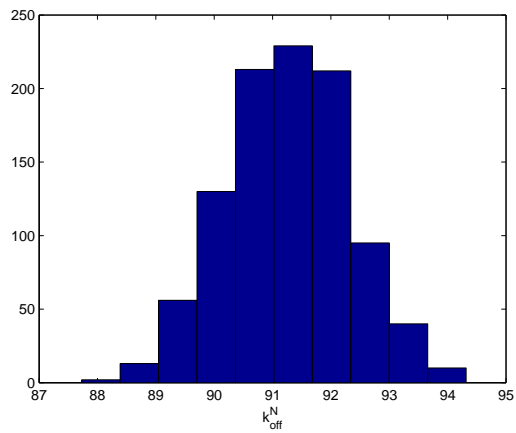}}
	\caption{Three parameters estimation ($k_I^+$, $k_I^-$ and $k_{off}^N$): Bootstrapping distribution for $k_{off}^N$. We used GLS and M=1000 runs.}	\label{fig:boot_1000_3p_koffn_koffn}
\end{figure}

\clearpage

\subsection{Estimation of 4 main parameters}

Following the sensitivity analysis detailed above, we tried to estimate a combination of the parameters $k_I^+, k_I^-, k_{on}^{N}, k_{off}^{N}$ for the parameter set with $\kappa_{\theta}=4$.\\
Parameters as follows were fixed from the original 9 parameter fit:
\begin{align*}
\begin{array}{|c|c |c |c |c|}
\hline
k_{on}^{min}& k_{on}^{max}& x_1 & x_2 & i_{max} \\
\hline
1684 & 1.5 \cdot 10^9 &  0.062 & 0.859 & 3.5\cdot 10^5\\
\hline
\end{array}
\end{align*}

We obtained the following result for the estimation of the four parameters using the data sets 1 to 4. In all of them, the condition number of the Fischer's information matrix $\kappa$ is too large to invert. This along with the sensitivity results above strongly suggests that the data sets do not contain sufficient information to estimate 4 or more parameters with any degree of certainty attached to the estimates.
\begin{align*}
\begin{array}{|c |c |c |c|c|c|c|}
\hline
 &k_I^+ & k_I^{-}& k_{on}^{N} & k_{off}^N & \sigma^2 & \kappa \\
\hline
DS 1 & 2.1431 & 12.4751  & 4616.962 & 108.259 &    8.7219\cdot 10^{-6} &    6.1226 \cdot 10^{19} \\
\hline
DS 2 & 2.7995 & 14.7630 & 4616.957& 108.4308  &  9.8694\cdot 10^{-6}  &  1.4442 \cdot 10^{19} \\
\hline
DS 3 &  2.180 & 15.757  & 4618.599& 41.369 &   6.4622\cdot 10^{-6} &   1.881\cdot 10^{17}\\
\hline
DS 4 &2.161 &10.9278 & 4617.3316 & 93.3265  &  6.374\cdot 10^{-6}&  2.144\cdot 10^{18}\\
\hline
\end{array}
\end{align*}




\section{Model Comparison Tests}\label{sec:chi2}

 A type of Residuals Sum of Squares (RSS) based model selection criterion \cite{BanksFitz,BHT2014, BT2009} can be used as a tool for model comparison for certain classes of models. In particular this is true for models such as those given in \cite{MComp} in which potentially extraneous mechanisms can be eliminated from the model by a simple restriction on the underlying parameter space while the form of the mathematical model remains unchanged. In other words, this methodology can be used to compare two nested mathematical models where the parameter set  $\Omega_{\theta}^{H}$ (this notation will be defined explicitly in Section \ref{sec.RSS.OLS} below) for the restricted model can be identified as a linearly restricted subset of the admissible parameter set $\Omega_{\theta}$ of the unrestricted model. Indeed, the RSS based model selection criterion is a useful tool to determine whether or not certain terms in the mathematical models are important in describing the given experimental data.

\subsection{Ordinary Least Squares}
\label{sec.RSS.OLS}
We now turn to the statistical model \eqref{absolute}, where the measurement errors are assumed to be independent and identically distributed with zero mean and constant variance $\sigma^{2}$. In addition, we assume that there exists $\theta_{0}$ such that the statistical model
\begin{equation}\label{eq.true.statistic.model}
Y_{j}=M(t_{j};\theta_{0})+\mathcal{E}_{j}, \quad j=1,2,\ldots,n.
\end{equation}
correctly describes the observation process. In other words, \eqref{eq.true.statistic.model} is the true model, and $\theta_{0}$ is the true value of the mathematical model parameter $\theta$.

With our assumption on measurement errors, the mathematical model parameter $\theta$ can be estimated by using the ordinary least squares method; that is, the ordinary least squares estimator of $\theta$ is obtained by solving
\begin{equation*}
\theta^{n} = \arg \min_{\theta \in\Omega_{\theta}} J_{}^{n}(\theta;\mathbf{Y}).
\end{equation*}
Here $\mathbf{Y} = (Y_{1}, Y_{2}, \ldots, Y_{n})^{T}$, and the cost function $J^{n}$ is defined as
\begin{equation*}
J^{n}(\theta;\mathbf{Y})  = \frac{1}{n}\sum_{k=1}^{n}\left(Y_{k} - M(t_{k};\theta)\right)^{2}.
\end{equation*}
The corresponding realization $\hat{\theta}^{n}$ of $\theta^{n}$ is obtained by solving
\begin{equation*}
\hat{\theta}_{}^{n} = \arg \min_{\theta \in\Omega_{\theta}} J^{n}(\theta;\mathbf{y}),
\end{equation*}
where $\mathbf{y}$ is a realization of $\mathbf{Y}$ (that is, $\mathbf{y} = (y_{1}, y_{2}, \ldots, y_{n})^{T}$).

As alluded to in the introduction, we might also consider a restricted version of the mathematical model in which the unknown true parameter is assumed to lie in a subset $\Omega_{\theta}^{H} \subset \Omega_{\theta}$ of the admissible parameter space.  We assume this restriction can be written as a linear constraint, $\mathcal{H}\theta_{0} = \mathbf{h}$, where $\mathcal{H} \in \mathbb{R}^{\kappa_{r} \times \kappa_{q}}$ is a matrix having rank $\kappa_{r}$ (that is, $\kappa_{r}$ is the number of constraints imposed),
and $\mathbf{h}$ is a known vector.  Thus the restricted parameter space is
\begin{equation*}
\Omega_{\theta}^{H} = \left\{\theta \in\Omega_{\theta} : \mathcal{H}\theta = \mathbf{h} \right\}.
\end{equation*}
Then the null and alternative hypotheses are
\begin{align*}
H_{0} & : \quad \theta_{0} \in \Omega_{\theta}^{H}\\
H_{A} & : \quad \theta_{0} \not\in \Omega_{\theta}^{H}.
\end{align*}
We may define the restricted parameter estimator as
\begin{equation*}
\theta^{n,H} = \arg \min_{\theta \in \Omega_{\theta}^{H}} J_{}^{n}(\theta;\mathbf{Y}),
\end{equation*}
and the corresponding realization is denoted by $\hat{\theta}_{}^{n,H}$.
Since $\Omega_{\theta}^{H} \subset \Omega_{\theta}$, it is clear that
$$J_{}^{n}(\hat{\theta}_{}^{n};\mathbf{y}) \leq J_{}^{n}(\hat{\theta}_{}^{n,H};\mathbf{y}).$$
This fact forms the basis for a model selection criterion based upon the residual sum of squares.
Using the standard assumptions  (given in detail in \cite{BHT2014}), one can establish asymptotic convergence result for the \textit{test statistics} (which is a function of observations and is used to determine whether or not the null hypothesis is rejected)
\begin{equation*}
U^{n} = \frac{n \left(J^{n}({\theta}^{n,H};\mathbf{Y}) - J^{n}({\theta}^{n};\mathbf{Y})\right)}{J^{n}(\theta^{n};\mathbf{Y})},
\end{equation*}
where the corresponding realization $\hat{U}_{n}$  is defined as
\begin{equation}\label{eq.RSS.OLS.statistic.realization}
\hat{U}_{}^{n} = \frac{n \left(J_{}^{n}(\hat{\theta}_{}^{n,H};\mathbf{y}) - J_{}^{n}(\hat{\theta}_{}^{n};\mathbf{y})\right)}{J_{}^{n}(\hat{\theta}_{}^{n};\mathbf{y})}.
\end{equation}
This asymptotic convergence result is summarized in the following theorem.
\begin{thm}\label{thm:model_comp2}
 Under assumptions detailed in \cite{BHT2014,BT2009} and assuming the null hypothesis $H_{0}$ is true, then $U_{}^{n}$ converges in distribution (as $n \to \infty$) to a random variable $U$ having a chi-square distribution\index{distribution!chi-square} with $\kappa_{r}$ degrees of freedom.
\end{thm}

The above theorem suggests that if the sample size $n$ is sufficiently large, then
$U_{}^{n}$ is approximately chi-square distributed with $\kappa_{r}$ degrees of freedom.
We use this fact to determine whether or not the null hypothesis $H_0$ is rejected.
To do that, we choose a \textit{significance level} $\alpha$ (usually chosen to be 0.05) and use
$\chi^2$ tables to obtain the corresponding \textit{threshold} value
$\tau$ so that $Prob(U>\tau) =\alpha$.
We next compute $\hat{U}_{}^{n}$ and compare it to $\tau$. If
$\hat{U}_{}^{n} >\tau$, then we {\em reject} the null hypothesis $H_0$ with confidence
level $ (1-\alpha)100\%$; otherwise, we do not reject.
 We emphasize that care should be taken in stating
conclusions: we either reject or do not reject $H_0$ at the
specified level of confidence. The table below illustrates the threshold values for $\chi^2(1)$ with the given significance level.
\begin{align*}
\begin{array}{|c| c| c|} \hline
            \alpha  & \tau & \text{confidence level}  \\ \hline
            .25   & 1.32   & 75\%  \\
            .1    & 2.71   & 90\%  \\
            .05   & 3.84   & 95\%  \\
            .01   & 6.63   & 99\%  \\
            .001  & 10.83  & 99.9\%\\ \hline
\end{array}
\label{chitable}
\end{align*}
 Similar tables can be found in any elementary statistics text or online or calculated by some software package such as Matlab, and is given here for illustrative purposes and also for use in the examples demonstrated below.

\subsection{Generalized Least Squares}
The model comparison results outlined can be extended to deal with generalized least squares problems in which measurement errors are independent with $\mathbb{E}(\mathcal{E}_{k}) = 0$  and  $Var(\mathcal{E}_{k}) = \sigma^{2}w^{2}(t_{k}, \hat{\theta})$, $k=1,2,\ldots,n$, where $w$ is some known real-valued function with $w(t,\hat{\theta})\neq 0$ for any $t$. This is achieved through rescaling the observations in accordance with their variance (as discussed in \cite{BHT2014}) so that the resulting (transformed) observations are identically distributed as well as independent.

\subsection{Results for PolyQ Aggregation Models}

We then carried out a series of model comparison tests (we again used DS 4) for nested models
to determine if an added parameter yields a statistically significantly improved model fit.
Our null hypothesis in each case was:  $H_0$: The restricted model is adequate (i.e., the fit-to-data is not significantly improved with the model containing the additional parameter as a parameter to be estimated).
We obtained the following results.

\begin{enumerate}
\item Model with estimation of $\{k_I^+, k_I^-\}$ vs. the model with estimation of $\{k_I^+, k_I^-, k_{off}^N \}$ :
We find with n=699,
$J_{n}(\hat{\theta}_{H}^{n};\mathbf{Y})= .0044192109$, $J_{n}(\hat{\theta}^{n};\mathbf{Y})= .0043709501$ and $\hat{U}_n=7.7178$. Thus we reject $H_0$ at a  99\% confidence level.

\item Model with estimation of $\{k_I^+, k_I^-\}$ vs. the model with estimation of $\{k_I^+, k_I^-, k_{on}^N \}$ :
We find $J_{n}(\hat{\theta}^{n};\mathbf{Y})=.0044192108$ with $\hat{U}_n=7.49\times 10^{-06}.$ Thus we don't reject $H_0$ at a 99\% confidence level.

\item Model with estimation of $\{k_I^+, k_I^-, k_{off}^N \}$ vs. the model with estimation of $\{k_I^+, k_I^-, k_{off}^N, k_{on}^N  \}$ :
To the order of computation we find no difference in the cost functions in this case and therefore we do not reject $H_0$ at a confidence level of 99\%.

\item Model with estimation of $\{k_I^+, k_I^-, k_{on}^N \}$ vs. the model with estimation of $\{k_I^+, k_I^-, k_{on}^N, k_{off}^N  \}$ :
We find $J_{n}(\hat{\theta}^{n};\mathbf{Y})=.0043709780$ with $\hat{U}_n=7.7133 $ and hence we reject $H_0$ with a confidence level of 99\%.
\end{enumerate}

From these and the preceding results we conclude the information content of the typical data set for the dynamics considered here will support at most 3 parameters estimated with reasonable confidence levels and these are the parameters $\{k_I^+, k_I^-, k_{off}^N \}$.

\section{Conclusions and Suggested Further Efforts}
For the efforts reported on above we make several conclusions.

For the majority of data sets, the GLS residual plots with $\gamma=0.6$ are random when fitted for data points $M(t_k) \geq 0.12.$ As conjectured earlier, this may be because the early formation of aggregates is somewhat stochastic in nature which is not well described by either the mathematical and/or statistical models. It appears that one needs special consideration of smaller polymer sizes. Indeed we suspect from additional discussions with our colleagues that perhaps the nucleation step might be dominated by a stochastic rather than deterministic process in the early stages (i.e., for small polymer sizes). This is a possible direction of further investigation.

Based on several different mathematical/statistical methodologies (sensitivities, asymptotic analysis, bootstrapping, model comparison tests), the data sets we considered do not contain sufficient information for the reliable estimation of all 9 parameters of interest. Indeed our findings suggest that at most 3 parameters can be reliably estimated with the data sets typical of those presented here, and that these parameters are  $\{k_I^+, k_I^-, k_{off}^N \}$. Recently related efforts \cite{OD} suggest that perhaps there are experimental design questions that could be addressed to collect data that might support the more sophisticated models derived in \cite{DouPri2012}, especially in order to investigate information coming from different initial concentrations. Indeed, we have considered here data sets related to experiments carried out with the same initial concentration. Adapting the previously used techniques to \emph{simultaneously} or \emph{successively} use all the information content in data sets carried out for different initial concentration is a challenging problem (see \cite{PP} for a discussion of the effect of initial concentration on nucleated polymerization).

Here we conclude that at most 3 parameters $\{k_I^+, k_I^-, k_{off}^N \}$ can be reliably estimated with the data sets investigated. The two first parameters determine the balance between the normal and abnormal protein concentrations and the third represents the stability of the nucleus against the degradation into monomeric entities. These three parameters are related to the early steps of the aggregation process, and thus we conclude that the model applied to these data sets does not provide any insight into the polymerization of larger polymers. Since this is the case, there is little motivation to modify the polymerization function depicted in Figure~\ref{fig:kon} until further data collection procedures are pursued.

\vspace{0.2cm}


\section*{Acknowledgements}
This research was supported in part (MD, CK) by the ERC Starting Grant SKIPPERAD, in part (HTB) by Grant Number NIAID R01AI071915-10 from the National Institute of Allergy and Infectious Diseases, and in part (HTB) by the
Air Force Office of Scientific Research under grant number AFOSR FA9550-12-1-0188.

\begin{thebibliography}{}

\end{thebibliography}


\begin{thebibliography}{9}

\bibitem{ABDR}
B.M. Adams, H.T. Banks, M.~Davidian, and E.S. Rosenberg,
\newblock Model fitting and prediction with {HIV} treatment interruption data,
\newblock Center for Research in Scientific Computation Technical Report CRSC-TR05-40, NC State Univ., October,
2005; {\em Bulletin of Math. Biology}, {\bf 69} (2007), 563--584.


\bibitem{OD} Kaska Adoteye, H.T. Banks and Kevin B. Flores, Optimal design of non-equilibrium experiments for genetic network interrogation, CRSC-TR14-12, N. C. State University, Raleigh, NC, September, 2014; {\em Applied Mathematics Letters}, {\bf 40} (2015), 84--89; DOI: 10.1016/j.aml.2014.09.013.

\bibitem{MComp} H.T. Banks, J.E. Banks, K. Link, J.A. Rosenheim, Chelsea Ross, and K.A. Tillman, Model comparison tests to determine data information content,
CRSC-TR14-13, N. C. State University, Raleigh, NC, October, 2014; {\em Applied Math Letters}, to appear.


 \bibitem{HIVUQ} H.T. Banks, R. Baraldi, K. Cross, K. Flores, C. McChesney, L. Poag, and E. Thorpe, Uncertainty quantification in modeling HIV viral mechanics,  CRSC-TR13-16, N. C. State University, Raleigh, NC, December, 2013;  {\em Math. Biosciences and Engr.}, submitted.


\bibitem{parameterselection} H.T. Banks, A. Cintron-Arias and F. Kappel, Parameter selection methods in inverse problem formulation, CRSC-TR10-03, N.C. State University,  February, 2010, Revised, November, 2010; in {\em Mathematical Modeling and Validation in Physiology: Application to the Cardiovascular and Respiratory Systems},(J. J. Batzel, M. Bachar, and F. Kappel, eds.), pp. 43 -- 73, Lecture Notes in Mathematics Vol. 2064, Springer-Verlag, Berlin 2013.

\bibitem{banks_modelling_2008}
     \newblock H. T. Banks, M. Davidian, S. Hu, G. M. Kepler, and E. S. Rosenberg, 
     \newblock {Modeling HIV immune response and validation with clinical data},
     \newblock {\em Journal of Biological Dynamics}, \textbf{2} (2008), 357--385.


\bibitem{BDK1} H.T. Banks, M. Doumic and C. Kruse, Efficient numerical schemes for Nucleation-Aggregation models: Early steps,   CRSC-TR14-01, N. C. State University, Raleigh, NC, March, 2014.

\bibitem{BanksFitz} H.T. Banks and B.G. Fitzpatrick, Statistical methods for model comparison in parameter estimation problems for distributed systems, {\it{Journal of Mathematical Biology}}, {\bf{28}} (1990), 501-527.




 \bibitem{BHT2014} H.T. Banks, S. Hu and W.C. Thompson, {\em Modeling and Inverse Problems in the Presence of Uncertainty}, Taylor/Francis-Chapman/Hall-CRC Press, Boca Raton, FL, 2014.


\bibitem{BT2009} H.T. Banks and H.T. Tran, {\em Mathematical and Experimental Modeling of Physical and Biological Processes}, CRC Press, Boca Raton, FL, 2009.


\bibitem{CL2010} V. Calvez and N. Lenuzza and M. Doumic and J.-P. Deslys and F. Mouthon and B. Perthame,
 {\em Prion dynamic with size dependency - strain phenomena}, {J. of Biol. Dyn.}, 4 (1), 28--42.

\bibitem{CR}
R.J. Carroll and D. Ruppert, {\em Transformation and Weighting in
Regression,} Chapman \& Hall, New York, 1988.

\bibitem{CWR}
R.J. Carroll, C.F.J. Wu and D. Ruppert,  The effect of
estimating weights in Weighted Least Squares, {\em J. Amer.
Statistical Assoc.}, {\bf 83} (1988), 1045--1054.




\bibitem{CGPV02} J.F. Collet, T. Goudon,F. Poupaud and A. Vasseur, The {B}ecker-{D}\"{o}ring system and its {L}ifshitz-{S}lyozov limit, {\em SIAM J. Appl. Math.}, {\bf 62} (2002), 1488--1500.

\bibitem{MDnotes}
M. Davidian, {\em Nonlinear Models for Univariate and Multivariate Response}, ST 762 Lecture Notes, Chapters 2, 3, 9 and 11, 2007; http://www4.stat.ncsu.edu/~davidian/courses.html


\bibitem{DG} M. Davidian and D.M. Giltinan, \emph{Nonlinear Models for Repeated Measurement Data}, Chapman and Hall, London, 2000.

\bibitem{DiCiEfron} 
T.J. DiCiccio and B. Efron, Bootstrap confidence intervals, {\em Statistical Science}, {\bf 11} (1995), 189--228.



\bibitem{DGL2009} M. Doumic, T. Goudon and T. Lepoutre, Scaling limit of a discrete prion dynamics model, {\em Commun. Math. Sci.}, {\bf 7} (2009), 839--865.


\bibitem{Efronbook} 
B. Efron, {\em The Jackknife, the Bootstrap and Other Resampling Plans}, CBMS 38, SIAM Publishing, Philadelphia, PA, 1982.


\bibitem{LM2002} P. Lauren\c{c}ot and S. Mischler, From the discrete to the continuous coagulation–fragmentation equations, {\em Proc. Royal Society of Edinburgh: Section A Mathematics}, {\bf 132} (2002), 1219--1248.


\bibitem{Lev02} R.J. LeVeque, {\em Finite-Volume Methods for Hyperbolic Problems}, Cambridge University Press, 2002.

\bibitem{PP} E.T. Powers and D.L. Powers, The kinetics of nucleated polymerizations at high concentrations: Amyloid fibril formation near and above the ``supercritical concentration'', {\em Biophysical J.}, {\bf 91} (2006), 122--132.


\bibitem{DouPri2012} S. Prigent, A. Ballesta, F. Charles, N. Lenuzza, P. Gabriel, L.M.  Tine, H. Rezaei and M. Doumic, An efficient kinetic model for assemblies of amyloid fibrils and its application to polyglutamine aggregation, {\em PLoS ONE}, {\bf 7} (2012), e43273; DOI:10.1371/journal.pone.0043273

\bibitem{Rezaei2007}
F. Eghiaian, T. Daubenfeld, Y. Quenet, M. van Audenhaege, A.P. Bouin, G. van der Rest, J. Grosclaude and H. Rezaei,
 {Diversity in {p}rion protein oligomerization pathways results from domain expansion as revealed by hydrogen/deuterium exchange and disulfide linkage}, {\em PNAS}, {bf 104} (18), 2007,  {7414--7419}.

\bibitem{Rubinow} S.I. Rubinow, {\em Introduction to Mathematical Biology},  John Wiley \& Sons, New York, 1975.



\bibitem{SeWi} 
G.A.F. Seber and C.J. Wild, {\em Nonlinear
Regression}, J. Wiley \& Sons, Hoboken, NJ, 2003. 





\bibitem{SHR2007} Wei-Feng Xue, S.W. Homans and S.E. Radford, Systematic analysis of nucleation-dependent polymerization reveals new insights into the mechanism of amyloid self-assembly, {\em Proc Natl Acad Sci U S A}, {\bf 105} (2008), 8926--8931.


\bibitem{Xue2009} W.-F. Xue, S.~W. Homans, and S.~E. Radford, Amyloid fibril length distribution quantified by atomic force
  microscopy single-particle image analysis,  {\em Protein Engineering, Design \& Selection:PEDS},
  {\bf 22} (2009),489--496.

\bibitem{Xue2013} W.-F. Xue and S.~E. Radford, An imaging and systems modeling approach to fibril breakage enables
  prediction of amyloid behavior,
{\em Biophysical {J}ournal}, {\bf 105} (2013), 2811--2819.


\end{thebibliography}

\end{document}